\documentstyle[11pt,amscd,amssymb,bbm,verbatim]{article}
 
\newtheorem{thm}{Theorem}[section]

\newtheorem{cor}{Corollary}[section]

\newtheorem{lemma}{Lemma}[section]

\newcounter{num}
\begin{document}

\font\bbb = msbm10
\def\Bbb#1{\hbox{\bbb #1}}
\font\lbbb = msbm7
\def\LBbb#1{\hbox{\lbbb #1}}
\def\squarebox#1{\hbox to #1{\hfill\vbox to #1{\vfill}}}
\newcommand{\qed}{\hspace*{\fill}
           \vbox{\hrule\hbox{\vrule\squarebox{.667em}\vrule}\hrule}\smallskip}

\newcommand{\beq}        {\begin{eqnarray}}
\newcommand{\eeq}        {\end{eqnarray}}
\newcommand{\be}         {\begin{equation}}  
\newcommand{\ee}         {\end{equation}}
\newcommand{\beqn}       {\begin{eqnarray*}}  
\newcommand{\eeqn}       {\end{eqnarray*}}

\def\ep{\varepsilon}
\def\R{\mathbb R}
\def\C{\mathbb C}
\def\a{{\alpha}}
\def\b{{\beta}}
\def\p{\partial}
\def\t{{\theta}}
\def\T{\mathbb T}
\def\C{{\mathbb C}}
\def\I{{\mathcal I}}
\def\V{{\rm{Var}}}
\def\CV{{\rm{Cov}}}
\def\l{{\ell}}
\def\ld{{\lambda}}
\def\G{{\Gamma}}
\def\d{{\delta}}

\title{Deviations from the Circular Law}

\author{B. Rider
\thanks{Department of Mathematics, Duke University,
Box 90320, Durham, NC 27708, USA. rider@math.duke.edu.
This work was supported in part 
by the grant NSF-9983320.}
}

\maketitle

\begin{abstract} 
Consider Ginibre's ensemble of $N \times N$ non-Hermitian  random
matrices in which all entries are independent  complex
Gaussians  of mean zero and variance $\frac{1}{N}$.  As $N \uparrow
\infty$ the normalized  counting measure of the eigenvalues converges
to the uniform measure on  the unit disk in the complex plane.  In this
note 
we describe fluctuations about this {\em Circular Law}. First we obtain
finite $N$ formulas for the covariance of certain linear statistics of the
eigenvalues.  Asymptotics of these objects coupled with a theorem
of Costin and Lebowitz then result in central limit theorems  
for a variety of these statistics.
\end{abstract}

\section{Introduction}
\setcounter{equation}{0}
\label{sec:int}

A fundamental non-Hermitian ensemble of Random Matrix Theory (RMT)
is that of $N \times N$ matrices with
independent complex Gaussian entries of mean zero and variance $\frac{1}{N}$.
This model is typically  attributed to Ginibre who derived (\cite{G})
the joint  distribution of eigenvalues  $\{z_k\}$ lying in the complex
plane $\C$. That object is given by
\begin{eqnarray}
\label{equ:def1}
{dP}_N(z_1,z_2, \dots, z_N) 
& = & \frac{1}{Z_N}
      \prod_{1 \le \l < k \le N} |z_{\l} - z_k|^2 
      \mu_N(dz_1) \cdots \mu_N(dz_N)
\end{eqnarray}
in which $\mu_N(dz) = e^{-N|z|^2} d \Re(z) d  \Im(z)$ and $Z_N$ is the
appropriate normalizer.  A straightforward Large Deviation analysis
will take you from (\ref{equ:def1}) to the {\em Circular Law}: for $N
\uparrow \infty$ the  measure $\frac{1}{N} \sum_{k=1}^N \delta_{z_k}$
tends weakly to the uniform measure on the disk $ |z| \le
1$.\footnote{Experts know this to hold under far more general
conditions on the underlying distribution, see \cite{B} for the best
result.}   In other words, the mean number of eigenvalues falling in
some subset of the disk is well approximated by $N \times$ the
normalized area.  The goal here is to describe the associated
fluctuations.

Our study is centered around linear eigenvalue statistics of the form
\[
    X(f) = \sum_{k=1}^N f(z_k),
\]
and the behavior of their covariances
\[
  {\CV}_N \Bigl(X(f), X(g)\Bigr)  = \int_{\C^N} X(f) X(g) d P_N -
   \int_{\C^N} X(f) d P_N  \int_{\C^N} X(g) d P_N
\] 
as $N \uparrow \infty$.  In particular, if the statistic is one of
either  the moduli, $f(z_k) = f(|z_k|)$, or the angles,  $f(z_k) =
f(\mbox{arg}(z_k))$, we obtain exact formulas for the covariance at finite
$N$ which are amenable to a precise asymptotic analysis.  Afterward,
noting that 
the Vandermonde component in (\ref{equ:def1}) identifies the present ensemble 
as an {\em determinantal point field}
\footnote{Also known as {\em fermionic fields}, these objects were introduced
by Macchi in \cite{Ma}; an excellent introduction to their applications
in RMT and beyond is contained in \cite{S4}.}, we may apply a basic
result of that theory due to Costin and Lebowitz along with
our asymptotics to conclude a central limit theorem for $X(f)$ in
several cases.

There is already a vast body of work on  fluctuations of linear
statistics in random matrix ensembles.  The majority of these concern
either Hermitian ensembles (for example, 
\cite{BS}, \cite{BE},  \cite{BW}, \cite{CL}, 
\cite{Gu}, \cite{I},
\cite{Jo1} \cite{SiS}, and \cite{S2}) 
or the
ensembles defined by Haar measure  on the classical compact groups
(\cite{DE}, \cite{KS}, \cite{Jo2}, \cite{S1}, and \cite{W} to name a few). 
Of course, works such as
\cite{S3} which focus on determinantal point 
fields take on a more general point of view. Nevertheless,   
in the  non-Hermitian setting it appears this type of 
fluctuation question is only directly considered in \cite{F}, and that 
more or less from a physics standpoint.  In fact, one  motivation for the present 
paper was the increased interest in the physics community in non-Hermitian 
ensembles (\cite{F} and \cite{FS} offer guides to that literature).  A 
second motivation
lies in the discovery that Ginibre's complex Gaussian ensemble
possesses a structure allowing for exact formulas that in turn 
may be studied via rather straightforward mathematical machinery 
(basically Laplace-type asymptotics).
Finally, it is expected that the results here should provide some indication as to 
the order and shape of fluctuations in
more general non-Hermitian ensembles.\footnote{Part of the sequel
carries over to Ginibre's Quaternion Gaussian ensemble, see Section 5. However, 
describing fluctuations in the real Gaussian ensemble will require 
a new approach due to the intricacies of that eigenvalue density, 
see \cite{LS} or \cite{E}.}

We begin with a  description of our covariance formulas. 
The case when $f$ and $g$ are functions of the spectral moduli alone 
is far the simpler of the two. The observation is that:

\begin{thm} 
\label{thm:cov1}
For statistics of the moduli, let  $f$ and $g$ be bounded over $\R_{+}$ then 
\beq
\label{equ:cov1} 
{\CV}_N \Bigl(X(f), X(g) \Bigr) =  \sum_{\l=1}^N 
 \Bigl\{ E \Bigl[ f
  \Bigl(\sqrt{\frac{1}{N} s_{\l} } \Bigr)  g \Bigl(\sqrt{ \frac{1}{N} s_{\l}}\Bigr) \Bigr] -  E
  \Bigl[  f \Bigl( \sqrt{ \frac{1}{N} s_{\l} } \Bigr) \Bigr]  E \Bigl[ g
  \Bigl(\sqrt{ \frac{1}{N} s_{\l} } \Bigr) \Bigr]  \Bigr\}
\eeq 
where $s_{\l} = \eta_1 + \cdots + \eta_{\l}$ for
$\{ \eta_{\l} \}$ a sequence of independent exponential random
variables of mean one.
\end{thm}

The angular statistics possess a more interesting structure, requiring
the setting down of some notation. The angles or arguments of
the $\{z_k\}$ are indexed to lie in the circle $\T = [-\pi,\pi]$, the
latter normalized throughout according to $\int_{\T} d \theta = 1$.
We recall the Dirichlet kernel
\[
   D_{\l}(\t) = \sum_{k = -\l}^{\l} e^{i k \t}  = \frac{\sin \Bigl(
           (\l +\frac{1}{2}) \t \Bigr)}{\sin \Bigl( \frac{1}{2} \t
           \Bigr)},
\]
which acts on functions $f$ of $\T$ by convolution:  $( D_{\l} \ast
f)(\t) = \int_{\T} D_{\l}(\t - \t^{\prime}) f(\t^{\prime}) \, d
\t^{\prime}$.  As is well known, this action produces an {\em
approximation of the identity}, that is,   $\lim_{k \uparrow \infty}  (D_k
\ast f)(\t)  = f(\t)$ at points of continuity of $f$.  Figuring more
prominently in what follows is a different kernel with that  property. We
define 
\be
\label{equ:ourker}
C_{\l}(\t)  =  2^{2 \l} \frac{ (\l !)^2}{ (2 \l)!}  \cos^{2 \l}(\t)
+   2^{2 \l+1} \frac{ (\G( \l + \frac{3}{2}))^2 }{  (2 \l + 1)!}
\cos^{2 \l + 1}(\t) 
\ee  
in which  $\G$ is the usual Gamma function. The result is:

\begin{thm} 
\label{thm:cov2}
For statistics of the spectral angles,  let $f,g \in L^2(\T)$.
We then have 
\beq
\label{equ:cov2}
  {\CV}_N
 \Bigl(X(f), X(g)\Bigr)  & = & N \Bigl( (f \ast
       \tilde{g})(0)  - \frac{1}{N} \sum_{\l=0}^{N-1}  (C_{\l} \ast f
       \ast \tilde{g})(0) \Bigr) \\  &  &  - \sum_{\l =
       \lceil \frac{N}{2} \rceil}^{N-1} \Bigl(  (C_{\l} -  D_{2N -2 \l - 2} \ast C_{\l}) \ast ( f
       \ast \tilde{g}) \Bigr)(0) \nonumber   
\eeq  in which
       $\tilde{g}(\t) = g(-\t)$. 
\end{thm}

Theorem \ref{thm:cov2} should be compared with the results of 
\cite{BDK} which discusses the allied question for Haar distributed
eigenvalues in the Unitary group $U(N)$.   As in \cite{BDK}, it is
striking that the covariance  depends on the test functions $f$ and
$g$ only  through $ f \ast {\tilde g}$.  Further, with the Cesar{\'o}
averages $ \frac{1}{N} \sum_{\l=0}^{N-1} C_{\l}(\theta)$ also forming 
approximation to the identity, the content of (\ref{equ:cov2}) is that
the large $N$ properties of the covariance are tied  to the error in
that approximation.  In \cite{BDK} fluctuations for  $U(N)$
are shown to be related in the same way to the error in Fej{\'e}r approximation.

Next we turn to the asymptotics.  
It is of interest to connect the growth of the covariance in $N$ to the smoothness
imposed on the underlying test functions  $f$ and $g$.
For example, the $O(1)$ fluctuations of
$X(f)$ for $f$ a bit better than twice differentiable plus a growth condition  in  the Hermitian
case (see \cite{Jo1}) and for $f \in H_2^{\frac{1}{2}}(\T)$ in $U(N)$ 
(see \cite{DE} or \cite{Jo2}) are well known manifestations of the rigidity of random 
matrix ensembles.  In line with those results, an immediate consequence of 
Theorem \ref{thm:cov1} is the following.

\begin{thm}  
\label{thm:r1}
Let $f(z) = f(|z|)$ and $g(z) = g(|z|)$ lie in $C^{2,\d}$ for $r = |z| \in (0,1+ \ep)$ with 
some  $\d > 0$ and $\ep > 0$ and otherwise be bounded.
Then, as $N \uparrow \infty$, 
${\CV}_N ( X(f), X(g) ) = \frac{1}{2}  
\int_0^1 f^{\prime}(r) g^{\prime}(r) r \, dr  + o(1) $. Moreover,
the centered (but unnormalized) random variable $X(f) - E[X(f)]$ converges 
to a mean zero Gaussian with the corresponding variance.
\end{thm}

Actually, the central limit theorem for $X(f(|z|))$ described in 
Theorem \ref{thm:r1} has already been discussed in \cite{F}. 
We in fact follow the ideas therein;  
a proof is included here for the
sake of completeness and rigor.  While \cite{F} employs the same product
structure behind Theorem \ref{thm:cov1}, the probabilistic
interpretation  which guides our results and allows for more to be
accomplished (see below) was apparently not noticed. 
Also in \cite{F}, the author conjectures (and provides heuristics
based on the associated log-gas) that the general (smooth) linear statistic 
should have order one  Gaussian fluctuations
in the large $N$ limit. Our next result shows that 
this is not the case.

\begin{thm} 
\label{thm:p1}
Let the function $f$ of the phase have Fourier coefficients
$\widehat{f}(k) = \int_{\T} e^{-i k \t}  f(\t) d \t$ satisfying
$\sum_{-\infty}^{\infty} k^{2} |\widehat{f}(k)|^2 < \infty$  and
likewise for $g$.  Then
\[
 {\CV}_N \Bigl( X(f), X(g) \Bigr) = \frac{\log N}{4} \int_{\T}
         f^{\prime}(\t) g^{\prime}(\t) \, d\t + o(\log N).
\] 
If a bit more smoothness on either test function is assumed, in
particular if $ \sum_{-\infty}^{\infty}$ $|k|^{2+\d}$
$|\widehat{f}(k)|^2 < \infty$ for some $\d >0$, the 
$o(\log N)$ error may be replaced by an $O(1)$.
\end{thm}

\noindent{Of} course, an $O(\log N)$ fluctuation is  rigid compared to the $O(N)$ characterizing 
an ensemble of independent particles.  Still, this marked difference between radial and 
angular statistics was not anticipated.

Changing focus, an important class of non-smooth test functions to consider
are indicators of given sets within the spectrum, thus 
providing an eigenvalue count. 
Setting this case of the {\em number statistic} apart, let us denote
 $\#{\Theta}{[\a,\b]} = X({\mathbbm{1}}_{ \{ {\mbox{arg}} (z_k) \in [\a, \b] \} } )$.
and  $ \#{\I}{[a,b]} = X({\mathbbm{1}}_{\{ |z_k| \in [a, b] \} } ) $ 
The next result concerns the growth of the variances 
of $\# \Theta[\a,\b]$ and $\# \I[a,b]$, with added attention to
the {\em mesoscopic scales}, that is,  when
$|\a-\b|$ (or $|a-b|$) $\downarrow 0$ as $N \uparrow \infty$. It is found
there is a break in the behavior depending on whether $\sqrt{N} |\a- \b|$
(or $\sqrt{N} |a-b|$) tends to $\infty$ or not.

For the angular number statistic the result is:

\begin{thm} 
\label{thm:num2} 
Note first by rotation invariance it is enough to consider 
the symmetric interval $[-\a/2, \a/2]$.  If $\a > 0$ remains
fixed while $N \uparrow \infty$, then
\be
\label{equ:v3}
    {\V}_N \Bigl( \# {\Theta} {[-\frac{\a}{2}, \frac{\a}{2}]}  \Bigr)  
  = \sqrt{N} \, \frac{1}{\pi^{3/2}} + O(\log N).
\ee
If instead $\a = \a(N) \downarrow 0$, $N \a(N) \uparrow \infty$, 
there are the following cases:
\be
\label{equ:v4}
  {\V}_N \Bigl( \# {\Theta} {[-\frac{\a}{2}, \frac{a}{2}]} \Bigr) =  \left\{ 
  \begin{array}{ll}  \sqrt{N} \frac{1}{\pi^{3/2}}(1+o(1)) & \mbox{ if }   
                   \lim_{N \uparrow \infty} \sqrt{N} \a =  \infty \\   
                     \sqrt{N} \frac{1}{\pi^{3/2}} I_{arg}(\b)  + O(\log N) 
      & \mbox{ if }  \lim_{N \uparrow \infty} \sqrt{N} \a = \b > 0 \\
                     N \a  (1+o(1)) & \mbox{ if }  \lim_{N \uparrow \infty} \sqrt{N} \a = 0,
  \end{array} \right.
\ee 
where, $I_{arg}(\b) \ge 0$  and is  such that 
$lim_{ \b \downarrow 0} I_{arg}(\b) = 0$ while $\lim_{\b \uparrow \infty} I_{arg}(\b) = 1$;
its definition may be found below in (\ref{equ:Iarg}). 
\end{thm}

The radial case has a related structure:

\begin{thm}  For fixed  $[a,b] \in (0,1)$ we have that, 
\label{thm:num1}
\be
\label{equ:v1}
  {\V}_N \Bigl( \#{\I} {[a,b]} \Bigr) = \sqrt{N} \, \Bigl(\frac{a+b}{ \sqrt{\pi}}\Bigr)  + O(1). 
\ee 
If rather $(b-a) \downarrow 0$ while $N (b-a) \uparrow \infty$, then
\be
\label{equ:v2}
  {\V}_N \Bigl( \# {\I} {[a,b]} \Bigr) =  \left\{ 
  \begin{array}{ll}  \sqrt{N}  \frac{a}{\sqrt{\pi}}(1+o(1)) & \mbox{ if }   
                   \lim_{N \uparrow \infty} \sqrt{N}( b - a) =  \infty \\   
                     \sqrt{N} \frac{a}{\sqrt{\pi}} 
                    I_{mod}(c) + O(1) & \mbox{ if }  \lim_{N \uparrow \infty} \sqrt{N}(b-a) = c > 0 \\
                    N (b^2 - a^2) (1+o(1)) & \mbox{ if }  \lim_{N \uparrow \infty} \sqrt{N}(b-a) = 0,
  \end{array} \right.
\ee 
in which, similar to the above,  $I_{mod} $ (defined in (\ref{equ:Imod})) is non-negative  and 
satisfies $ \lim_{c \downarrow 0} I_{mod}(c) = 0$, $\lim_{c \uparrow \infty} I_{mod}(c)= 1$.
\end{thm}

\noindent{\bf Remark} Note that in the radial case the variance asymptotics
are modulated by the limiting radial distribution, the corresponding angular distribution
being uniform.

\bigskip

We may now invoke the Costin-Lebowitz Theorem (see Section 5).  
First proved in \cite{CL} for the sine kernel and extended by A. Soshnikov in
\cite{S2} and \cite{S3}, in the present  setting the theorem implies 
that: with $N \a \uparrow \infty$ (respectively $N (b-a) \uparrow \infty$)
\be
\label{equ:rat}
\frac{\# \Theta [-\a,\a] - E_N [ \# \Theta [-\a, \a]]}{\sqrt{{\V}_N( \# \Theta [-\a,\a])}} 
\ \ \ \left(\mbox{respectively  }
\frac{\# \I [a,b] - E_N [ \# \I [a,b]]}{\sqrt{{\V}_N( \# \I [a,b])}} \right)
\ee
tends to a Gaussian random variable of mean zero and variance one.  
Since we may also compute $E_N [ \# \Theta [-\a, \a]] =  N(2 \a)$
and $ E_N [ \# \I [a,b]] = N (b-a) + O(1)$, an interesting observation
regarding the mesoscopic scales is that when $
\lim_{N \uparrow \infty} \sqrt{N} \a < \infty$
(or  $\lim_{N \uparrow \infty} \sqrt{N} (b-a) < \infty$) we have
that ${\V}_N  [ \# \Theta [-\a, \a]] = O ( E_N [ \# \Theta [-\a, \a]] )$
(or  $ {\V}_N  [ \# \I [a, b]] = O ( E_N [ \# \I [a, b]] ) $) as in the
case of independent particles.  Contrariwise, if 
$ \lim_{N \uparrow \infty} \sqrt{N} \a = \infty$
(and  $\lim_{N \uparrow \infty} \sqrt{N} (b-a) = \infty$)
the slow growth of the variance compared to the mean indicates 
rigidity.  These remarks are suggestive of a possible asymptotic
independence taking place on the small scales.  In this context
we add that if one considers the outlying or edge eigenvalues 
(which may be considered a small scale) it is found that
$ {\V}_N [ \# \I [1,\infty] ] = O (E_N [ \# \I [1, \infty]]) = O(\sqrt{N})$.
Further evidence of an asymptotic independence of the
edge is contained in  \cite{R1} and \cite{R2} which show the largest
eigenvalues in absolute value respond to a limit theorem shared by 
independent sequences (see also Section 5).

Finally, consider the limiting Gaussian field $ \Theta_{\infty} [\a, \b]$
on $\T$ (or $ \I_{\infty}[a,b]$ on  $[0,1]$) resulting from the 
properly normalized statistics (\ref{equ:rat}) for fixed arguments
$(\a, \b)$ (or $(a, b)$) as $N \uparrow \infty$.
The corresponding correlations are anticipated to be such that, taking 
the first case, 
$ \Theta_{\infty} [\a, \b]$ and $ \Theta_{\infty} [\d, \gamma]$
are independent if either $[\a, \b]$ and $[\d, \gamma]$ are disjoint
or one is properly contained in the other; those intervals can positively 
or negatively correlated when they share a single endpoint. 
This rather odd structure was first observed by Wieand (see \cite{W}) for $U(N)$ and 
has since been rediscovered in this and other ensembles by a variety of methods.
If we wanted only to obtain this limit, an extension of Costin-Lebowitz 
could be employed (see \cite{S3}).
However, another advantage of the formulas obtained
in Theorems \ref{thm:cov1} and \ref{thm:cov2}  
is that we can directly describe the correlations at all large but finite
values of $N$.  Our last result is then as follows: note the 
shape of the second order terms which speak to the differing degree 
of rigidity in the two cases,

\begin{thm}  
\label{thm:lcov}
The covariance of   $\# \Theta [\a,\b]$ and $\# \Theta [\delta,\gamma]$ 
satisfies  
\[ 
{\CV}_N  \left(
\frac{ \# \Theta [\a,\b]}{(N \pi^3 )^{1/4}} ,
\frac{ \# \Theta [\delta,\gamma]}{(N \pi^3 )^{1/4}} 
\right) = \left\{
\begin{array}{ll} \pm O(\frac{1}{\sqrt{N}}) & \mbox{if } [\a,\b]
\subset  [\delta,\gamma]  \mbox{ or } [\a,\b]  \cap [\delta,\gamma] =
\varnothing,  \\ + \frac{1}{2} + O(\frac{\log N}{\sqrt{N}} )     &
\mbox{if } \a = \delta \, (\b \neq \gamma) \mbox{ or } \b = \gamma
                        \, (\a \neq \d),  \\ 
- \frac{1}{2} +
O(\frac{\log N}{\sqrt{N}} )     & \mbox{if } \b = \delta \, (\a \neq \gamma) \mbox{ or }
\a = \gamma \, (\b \neq \d) .   \end{array} \right.  
\] 
For  $\# \I {[a, b]} $ and $\# \I{[c,d]}$ the result is
\[
{\CV}_N  \left(
\frac{\# \I [a,b]}{(N \pi )^{1/4}} ,
\frac{\# \I [c,d]}{(N \pi)^{1/4}} \right) = \left\{
\begin{array}{ll} \pm O( e^{-c_1 N} ) & \mbox{if } [a,b]
\subset  [c,d]  \mbox{ or } [a,b]  \cap [c,d] =
\varnothing,  \\ + \frac{a}{2} + O(\frac{1}{\sqrt{N}} )     &
\mbox{if } a = c \, (b \neq d) \mbox{ or } b = d \, (a \neq c),  \\ - \frac{b}{2} +
O(\frac{1}{\sqrt{N}} )     & \mbox{if } b = c  \, (a \neq d) \mbox{ or }
a = d  \, ( b \neq c) \end{array} \right.  
\] 
where $c_1 > 0$ and depends on $a,b,c,d$.
\end{thm}

As for the remainder of this note: in the next section we
derive our basic formulas for the covariance  (Theorems \ref{thm:cov1}
and \ref{thm:cov2}), the verification of  the asymptotic statements
are found in Section 3 (angular case)  and Section 4 (radial case).
Section 5 extends the present results in the radial setting 
to Ginibre's Gaussian quaternion ensemble.
For the reader's convenience Section 6 serves as a brief appendix 
on the Costin-Lebowitz theorem.

\bigskip 

\noindent{\bf Remark}
In a recent extension of Costin-Lebowitz,
\cite{S2} shows that (modulo technicalities) linear statistics in general determinantal
point fields satisfy a central limit
theorem  as long as the variance  grows faster than a small positive
power of the expectation.  The estimates below then imply a central
limit theorem at rate $N^{1/4}$ for $X(f)$ (either in the angular or
radial case) when  $f$ has anywhere a jump
discontinuity. Unfortunately, the  $O(\log N)$ growth of the variance
in Theorem {\ref{thm:p1}}  is not fast
enough from this point of view.  A proof of the central limit theorem
in the case of smooth angular statistics  
will appear elsewhere.

\section{Covariance Formulas}
\setcounter{num}{2} \setcounter{equation}{0}
\label{sec:vcv}

Appropriate row and column operations in (\ref{equ:def1})  allow
the eigenvalue density to be re-expressed as in 
\be
\label{equ:def2}
d {P}_N(z_1,z_2, \dots z_N)  =   \frac{1}{N!} \det \Bigl[ K_N (z_{\l},
 \bar{z_k})  \Bigr]_{1 \le \l,k \le N} \mu_N(dz_1) \cdots \mu_N(dz_N)
\ee 
with the Hermitian kernel 
\[
 K_N(z, \bar{w})  = 
 \sum_{\l = 0}^{N-1} \frac{N^{\l+1}}{\pi \l !}  z^{\l} {\bar{w}}^{\l}, 
\]
see \cite{M}.  As an operator on $L^2(\C, \mu_N(dz))$,
 $K_N(z, \bar{w})$ projects onto the span of
the first $N$ polynomials orthogonal with respect to that weight
(those being just the monomials  $1, z, z^2$ etc.).  This fact explains the rules 
$\int_{\C} K_N(z, \bar{z}) \mu_N(dz) = N$
and $\int_{\C} K_N(z, \bar{v}) K_N( v, \bar{w}) \mu_N(dv) =  K_N(z,
\bar{w})$, from which it follows that the marginal densities of 
$P_N$ are given by: 
\beq
\label{equ:cor}
\lefteqn{P_N^k(z_1, \dots, z_k)} \\ 
& := &  
 \Bigl\{ \int_{\C^{N-k}}
    \frac{1}{N!} \det \Bigl[ K_N(z_i, \bar{z_j}) \Bigr]_{1\le i,j \le N}
 \mu_N(dz_{k+1}) \cdots \mu_N(d z_N) \Bigr\} \mu_N(d z_1) \cdots \mu_N(dz_k)  \nonumber \\    & = &
 \frac{(N-k)!}{ N! }  \det \Bigl[ K_N(z_i, \bar{z_j}) \Bigr]_{1 \le
 i,j \le k} \mu_N(dz_1) \cdots \mu_N(dz_k).  \nonumber 
\eeq 
(normalized differently these objects are also known as the correlation 
functions).  Applying the integrating-out rules behind (\ref{equ:cor}) along 
with the symmetries of the integrands below we find  
\be
 E_N \Bigl[ X(f) \Bigr]  = N \int_{\C^N} f(z_1) dP_N = \int_{\C} f(z)
 K_N(z,\bar{z})  \mu_N(dz)
\ee
and 
\beq
\label{equ:covfin}
{\CV}_N \Bigl( X(f), X(g) \Bigr)   &  = & 
 \int_{\C} f(z)g(z) K_N(z,\bar{z}) \mu_N(dz)  \\
 & &      - \int_{\C} \int_{\C}  f(z) g(w) K_N(z,\bar{w}) K_N(w, \bar{z})
 \mu_N(dz) \mu_N(dw).  \nonumber
\eeq 
With that our expressions for the covariance may be established.

\bigskip

\noindent{\bf Proof of Theorem \ref{thm:cov1}}  In polar coordinates,
$z  = (r, \a), w = (s, \b)$, you find that
\[
  \int_{\C} f(|z|) g(|z|) K_N(z, \bar{z}) \mu_N(dz)  =  \sum_{\l
  =0}^{N-1} N^{\l+1} \frac{2}{\l !} \int_0^{\infty} f(r) g(r) r^{2 \l + 1}  e^{-N
  r^2} \, dr,
\] 
while in the second integral of  (\ref{equ:covfin}) only the diagonal
term survives: 
\beqn 
\lefteqn{\int_{\C} \int_{\C}  f(|z|)g(|w|)
|K_N(z, \bar{w})|^2 \mu_N(dz) \mu_N(dw)} \\ & = & \sum_{\l = 0}^{N-1}
\Bigl[ 4 \frac{N^{2 \l +2}}{(\l !)^2}  
     \int_0^{\infty} \int_0^{\infty} f(r) g(s) r^{2 \l + 1} s^{2 \l + 1}  e^{-N r^2}
         e^{-N s^2} \, dr \, ds \Bigr].   
\eeqn 
The connection with
(\ref{equ:cov1}) is then made after a change of variables: the
distribution  $ p(dr) := (N^{\l+1}/\l!) r^{\l} e^{-N r} \, dr $  being
recognized as that for a sum of $(\l+1)$ independent exponential random
variables, each of mean $\frac{1}{N}$.  The proof is finished.

\bigskip 

\noindent{\bf Proof of Theorem \ref{thm:cov2}}  From here on we denote
$\phi(\t) = (f \ast {\tilde g})(\t)$.  Also, in the present derivation
we consider $N$ to be even just to have things set;  the formulas are easily 
adjusted for odd values of $N$.

To arrive at (\ref{equ:cov2}), it
is now the radial component which is integrated out. Doing so in 
the first term
of (\ref{equ:covfin}) is immediate:
\[
  \int_{\C} (fg)(\mbox{arg}(z)) K_N(z,\bar{z}) \mu_N(dz) = \int_{\T} f(\t)
  g(\t) d \t  \int_0^{\infty} K_N(r,r) r e^{-Nr^2} dr  = N \phi(0).
\]
While for the second term we have
\beqn \lefteqn{ \int_{\C} \int_{\C} f(\mbox{arg}(z))
g(\mbox{arg}(w)) |K_N(z,\bar{w})|^2  \mu_N(dz) \mu_N(dw)} \\ & = &
\int_{\T} \int_{\T} f(\a) g(\b) \Bigl[  \sum_{0 \le k,\l \le N-1}
e^{i \a(k-\l)  }  e^{-i \b (k - \l)  }
\frac{ ( \int_0^{\infty} r^{k+\l+1} e^{-Nr^2} dr )^2}{ k !
\l!}  N^{k+\l+2}  \Bigr] d\a d \b.  
\eeqn 
With $\widehat{\phi}(k) =
\int_{\T} e^{i k \t} \phi(\t) d \t = \widehat{f}(k)
\widehat{g}(-k)$ the $k$-th Fourier coefficient of $\phi$, this produces
\beq
\label{equ:firstsum} 
  {\CV}_N \Bigl(X(f), X(g) \Bigr)  & = &      N \phi(0) - \sum_{0
\le k, \l \le N-1}  \Bigl[ \Gamma \Bigl( \frac{k+\l}{2} + 1 \Bigr) \Bigr]^2
\frac{\widehat{\phi}(k-\l)}{ k! \l !}  \\ & := &  N \phi(0) -
S_N^{0}(\phi). \nonumber   
\eeq 
as a preliminary form of our covariance formula. The obvious next 
step is make the substitution $ (k
+ \l, k-\l)  = (2n, 2m)$ for $k$ and $\l$ of the  same parity and
$(k+\l, k-\l) = (2n + 1, 2m-1)$ otherwise.  The result of that move is
\be
\label{equ:secsum}
   S_N^{0}(\phi) = \sum_{(m,n) \in {{\cal D}}_N}  \frac{(n!)^2
                   \widehat{\phi}(2m)}{(n-m)!(n+m)!}  + \sum_{(m,n) \in
                   {{\cal D}_{N}^{\prime}}} \frac{( \Gamma(n +\frac{3}{2}) \, )^2
                   \widehat{\phi}{(2m-1)}}{(n+m)! (n-m+1)!} 
\ee 
in which the sets ${\cal D}_N$ and ${\cal D}_N^{\prime}$
in the $(m,n)$-plane are described as follows. ${\cal D}_N$
consists of two triangles: the lower defined by the points 
$(0,0)$, $(N/2 -1, N/2 - 1)$ and $( -N/2+1, N/2-1)$, and the upper
defined by $(0, N-1)$, $(N/2 -1, N/2)$ and $(-N/2+1, N/2)$. The
set ${\cal D}_N^{\prime} = {\cal D}_N \cup {\cal D}_N^{\prime \prime}$
where the latter consists of all integer points lying on the
lines between $(1,0)$ and $(N/2+1, N/2 -1)$ and between
$(N/2 +1, N/2)$ and $(1, N-1)$.

To evaluate (\ref{equ:secsum}), first consider the sum over the lower
half of  ${\cal D}_N$. In particular start with a  fixed  $n \le N/ 2 - 1 $ 
so that the corresponding 
sum over $m$ ranges fully between $-n$ and $+n$. This
object may be computed by the following observation.
Denote by $\zeta_1, \zeta_2,\dots $ a sequence of independent $\pm1$  Bernoulli
random variables with mean $E$, and notice that:
with $c(n) = 2^{2n} {{2 n} \choose n}^{-1}$, 
\beq
\label{equ:even}
\lefteqn{\sum_{m=-n}^{n} \frac{(n!)^2 \,  \widehat{\phi}(2m)}{(n-m)!(n+m)!}
 =    c(n) 
\sum_{m=-n}^n  \widehat{\phi}(2m) \,  { {2 n} \choose n + m } 2^{-2n} } \\ 
& = & c(n)  E \Bigl[ \widehat{\phi}(\zeta_1 + \cdots +
 \zeta_{2n}) \Bigr]   
 =   c(n)      \int_{\T}
        \Bigl( E \Bigl[ e^{\sqrt{-1} \zeta_1 \t}  \Bigr] \Bigr)^{2n}
        \phi(\t) d \t  =   c(n)  \int_{\T} \phi(\t) \cos^{2n}(\t) \, d \t.  \nonumber 
\eeq 
In the same way we have  
\be
\label{equ:odd}
\sum_{m=-n}^{n+1} \frac{ (\G(n + 3/2) )^2 \widehat{\phi}(2m+1)}{(n-m)!(n+m+1)!}  
 =  2^{2n+1}  \frac{(\G(n+ 3/2))^2}{(2n+1)!}  \int_{\T} \phi(\t) \cos^{2n+1}(\t) \, d \t
\ee 
for the typical term in lower half of the ${\cal D}_N^{\prime}$ sum. Pairing
the final expressions in (\ref{equ:even}) and (\ref{equ:odd}) in $n$ brings out 
the sum of $(C_n
\ast \phi)(0)$ for $n$ from $ 0$ to $ N/2 - 1$.

The stated form of the covariance is then obtained by extending the
summation over $m$  in the upper reaches of the $n$ variables (extending
to a ``big triangle''), in order to produce a summation of $C_{n}$ over all
$n \le N-1$.  That is, we write
\[
   S_N^{0}(\phi)  := S_N(\phi) + S_N^{\prime}(\phi)
\]
with $ S_N(\phi) = \sum_{n = 0}^{N-1} (C_{n} \ast \phi)(0)$ and
$S_N^{\prime }(\phi)$ the error involved in throwing too much into
the mix.  To complete the picture,
we spell out the contribution to $S_N^{ \prime}(\phi)$  
connected to the equal-parity, or ${\cal D}_N$, sum.  That reads:
with $p_{2n}(m) = { {2n} \choose{ n + m}} 2^{-2n}$ and $c(n)$
as before, 
\beqn
\lefteqn{\sum_{n = \frac{N}{2}}^{N-1} c(n) \sum_{m = -n}^n  
          \widehat{\phi}( 2 m) \Bigl[ 1 - {\mathbbm 1}_{\{ |m| \le N - 1 - n \}} \Bigr] \, p_{2n}(m) } \\
& =  &   \sum_{n = \frac{N}{2}}^{N-1} c(n) \, \sum_{m= -n}^n 
                        \Bigl[ \Bigr(   \phi  
                          - ( D_{ 2 N -2n - 2  } \ast \phi ) \Bigl) {\widehat{ \phantom{h}}}  
             \hspace{-.1cm} ( 2m ) \Bigr] \, p_{2n}(m)  \\
& =   &  \sum_{n = \frac{N}{2}}^{N-1}  2^{2n}
{ {2n} \choose n}^{-1} \hspace{-.2cm}  \int_{\T} \cos^{2n}(\t) \Bigl[  \phi(\t) - ( D_{ 2 N -2 n - 2} \ast \phi )(\t) 
                   \Bigr] \, d \t.
\eeqn
The sum corresponding to the ${\cal D}_N^{\prime}$ part of our expression
is similar; that its outcome is as advertised should be clear. The proof is 
finished.

\section{Angular Statistics}
\setcounter{num}{3} \setcounter{equation}{0}
\label{sec:angs}

That the kernels $C_{\l}$ form an approximation of the identity is
seen from: (1) the evaluation $\int_{\T} C_{\l}(\t) \, d \t = 1$ 
(based on 
$\int_{\T} \cos^{2\l} (\t) \, d\t = 2^{-2\l}{ {2 \l} \choose \l}$ 
and $\int_{\T} \cos^{2\l +1 } (\t) \, d \t = 0$), and (2), the basic estimate 
\be
\label{equ:dummy}
  \int_{-\ep}^{\ep}  C_{\l}(\t) d \t 
 =  \sqrt{\frac{\l}{\pi}} \int_{-\ep}^{\ep} e^{-\l \t^2} d \t + o(1) = 1 + o(1)
\ee
valid for any small $\ep >0$.  The proofs of 
the main theorems hinge on the justification of exactly this sort of
standard Laplace asymptotics, though it is necessary to extract the
precise form of the error terms.  For the applications in mind, $C_{\l}$
acts on one of two classes of test functions. We restrict our attention
to those cases:

\begin{lemma} 
\label{lem:error}
Let $h(\t)$ be twice continuously differentiable on $\T$. Then, as $\l \uparrow \infty$,
\be
\label{equ:er1}
   (C_{\l} \ast h)(0) = \int_{\T} h(\t) C_{\l}(\t) \, d \t = h(0) + \frac{ h^{\prime \prime}(0)}{4\l}  + 
            \frac{\nu(\l, h)}{\l} + O \Bigl(\frac{1}{\l^{3/2}} \Bigr).  
\ee 
In general, $\nu(\l,h) = o(1)$ for large values of $\l$. This may be improved to
$\nu(\l,h) = O(\l^{-\d/2})$ if $h^{\prime \prime}(\t)$ happens to be
H\"older continuous of order $\d >0 $
in a neighborhood of the origin.

Let instead $h(\t)$ be a continuous, piecewise-linear function of $\t \in \T$. In
that case,
\be
\label{equ:er2}
  (C_{\l} \ast h)(0)  = h(0) + \frac{h_{+}^{\prime}(0)}{2
   \sqrt{\pi \l}}  - \frac{h_{-}^{\prime}(0)}{2 \sqrt{\pi \l}}  + O
   \Bigl( \frac{1}{\l^{3/2}} \Bigr)
\ee
is the appropriate estimate for $\l \uparrow \infty$. Here $h_{+}^{\prime}$
and $h_{-}^{\prime}$ indicate right and left derivatives. 
\end{lemma}

A comment is perhaps in order at this point as to the comparison 
between the kernel arising and being analyzed here and the better
known Dirichlet kernel.  The content on the second equality in 
(\ref{equ:dummy}) is that  $C_{\l}(\t)$ concentrates its mass in a neighborhood
of width $O(\l^{-1/2})$.  The concentration of the Dirichlet is sharper,
taking place in a neighborhood of order $\l^{-1}$. Returning to the discussion 
after Theorem \ref{thm:cov2} this explains why the
fluctuations for the angles in Ginibre's ensemble may 
be expected to be larger than those for $U(N)$.

\bigskip

\noindent{\bf Proof of Lemma \ref{lem:error}} For $h$ in either class, it may be assumed from
the start that $h(0) = 0$ ($\int_{\T} C_{\l}(\t) \, d \t = 1$). 
Also, throughout the proof we will make use of the simple
estimates
\be 
\label{equ:stfacts}
  2^{2 \l} { {2 \l} \choose \l }^{-1} = {\sqrt{\pi \l} } 
           \Bigl( 1 + \frac{1}{8 \l} + O \Bigl( \frac{1}{\l^2} \Bigr) \Bigr)
\ \mbox{ and  } \ 
  \frac{\Gamma(\l + \frac{1}{2})}{\l !} = \frac{1}{\sqrt{\l}} 
                  \Bigl( 1 - \frac{1}{8 \l} + O \Bigl( \frac{1}{\l^2} \Bigr) \Bigr),
\ee
which follow from  Stirling's 
approximation  in the form $  \Gamma(\l) = {\l}^{\l- \frac{1}{2}} e^{-\l} \sqrt{2 \pi} (   1 +
\frac{1}{12} {\l}^{-1} + O(\l^{-2}) )$.

Making the abbreviation 
\beqn
C_{\l}(\t)  & = &  \Bigl[ 2^{2 \l} \frac{(\l !)^2}{ ({2 \l})!}  \Bigr] \cos^{2\l} (\t) + 
              \Bigl[ 2^{2 \l+1}  \frac{\Gamma^2(\l + \frac{3}{2})}{(2 \l + 1) !} \Bigr] 
               \cos^{2 \l +1}(\t) \\ 
& := & a(\l) \cos^{2 \l} (\t) + b(\l) \cos^{2 \l + 1}(\t),
\eeqn 
the rule $\Gamma(k) = (k-1) \Gamma(k-1)$ provides the identity
$b(\l) = a(\l) \times (\l + \frac{1}{2}) \times (\Gamma(\l + \frac{1}{2}) / \l !)^2 $,
and from (\ref{equ:stfacts}) we note that
$a(\l) \simeq b(\l) \simeq \sqrt{\pi \l}$. The basic behavior of $C_{\l}$ 
for $\l \uparrow \infty$ may then be described as a 
sum of two (approximate) point masses of weight $\simeq 1/2$ at $\t = 0$ and the difference of two such 
masses at the common point $\t = \pm \pi$.  
Our first step is to dispense of the cancellation which must take place at the latter point.

Consider in particular 
the integral $\int_{\pi/2}^{\pi} h(\t) C_{\l}(\t) \, d \t$,  with that over $\t \in [-\pi, -\pi/2]$
being treated similarly.  Denoting  $g(\t) = h(\t) - h(\pi)$ we write
\beq
\label{equ:cancel}
  \int_{\pi/2}^{\pi} C_{\l}(\t) h(\t) \, d \t & =  &  h(\pi) \int_{\pi/2}^{\pi} C_{\l}(\t) d \t 
   +  (a(\l) - b(\l)) \int_{\pi/2}^{\pi} \cos^{2 \l}(\t) g(\t) \, d \t  \\
  &  &  \ +   b(\l) \int_{\pi/2}^{\pi} \cos^{2 \l}(\t) g(\t) ( 1 + \cos(\t)) \, d \t.  \nonumber
\eeq
This object is to be dominated by a constant multiple of $\l^{-3/2}$.  The first term on the
right hand side responds to an exact calculation plus by an appeal to (\ref{equ:stfacts}):
\beqn
\int_{\pi/2}^{\pi} C_{\l} (\t) d \t 
& = &  a(\l) \int_0^{\pi/2} \cos^{2 \l} (\t) d \t 
                    - b(\l) 
                               \int_0^{\pi/2} \cos^{2 \l +1}(\t) d \t  \\
& =   &  \frac{\pi}{2} - \frac{1}{2} 
                  \Bigl(  2^{2 \l} { {2 \l} \choose \l }^{-1} 
                          \times \frac{ \Gamma(\l + \frac{1}{2})}{\l !}  \Bigr)^{2}
           = O \Bigl( \frac{1}{\l^2} \Bigr). 
\eeqn
As for terms two and three, 
our regularity assumptions on $h$ imply that $|g(\t)| \le c_1 |\pi - \t|$
and so $| g(\t) (1 + \cos(\t) )| \le c_2 |\pi - \t|^3$ for $\t \in [\pi -\alpha, \pi]$ with a 
small $\alpha > 0$. Also on that same range we have the  bound $\cos^{2\l}(\t) \le e^{-c_3 \l (\t - \pi)^2}$
for $c_3 = c_3(\alpha) > 0$, while a fixed distance away from $\t = \pi$ (or $\t = 0$) 
$\cos^{2 \l}(\t)$ is exponentially small. It follows that
\[
 \Bigl| \int_{\pi/2}^{\pi} \cos^{2 \l}(\t) g(\t)  d \t \Bigr| \le c_4 \frac{1}{\l} 
\ \mbox{ and } \
  \Bigl| \int_{\pi/2}^{\pi} \cos^{2 \l}(\t) g(\t) ( 1 + \cos(\t)) \, d \t \Bigr| \le c_5 \frac{1}{\l^2}.
\]
Using again (\ref{equ:stfacts}) shows that
\[
a(\l) - b(\l)  =   2^{2 \l} { {2 \l} \choose \l }^{-1}  
                  \Bigl[ 1 - (\l + \frac{1}{2}) \Bigl( \frac{ \Gamma(\l + \frac{1}{2})}{\l !} \Bigr)^2
                     \Bigr] = O\Bigl( \frac{1}{\sqrt{\l}} \Bigr)  
\]
and completes our consideration of (\ref{equ:cancel}).

Turning our attention to the contribution to the integral from the
vicinity of $\t = 0$, take first the case that $h \in C^2$ and
use Taylor's theorem with remainder to write
\be
\label{equ:taylor}
   h(\t) = h^{\prime}(0) \t  + \frac{1}{2}
                      h^{\prime \prime}(0) \t^2 + \int_0^{\t}
                      \Bigl( h^{\prime \prime} (\a) - h^{\prime
                      \prime}(0) \Bigr) (\t -\a) d \a,
\ee
(recall $h(0) = 0$).  Since $C_{\l}(\t)$ is even there is no contribution
from the first term on the right hand side: $\int_{-\pi/2}^{\pi/2} C_{\l}(\t) \t \, d \t = 0$. 
Continuing to the quadratic term, one may again cut down the range of integration.
Needing to be more precise this time around  we restrict to $| \t | \le \l^{-1/4}$,
the ensuing error still
exponentially small. Further, on that range 
it holds that $| 1 - \cos(\t)| \le
\t^2$ and  $ |1- e^{\l \t^2} \cos^{2 \l}(\t) | \le 4 \l \t^4$.  Using
these estimates in $\int_{-\pi/2}^{\pi/2} C_{\l}(\t) \t^2  \, d \t = 
\int_{-\l^{-1/4}}^{\l^{-1/4}} C_{\l}(\t) \t^2  d \t + O(e^{-c_6 \l}) $    
and then restoring the
limits of integration over the whole line, we have   
\beqn 
\lefteqn{ 
\int_{-\pi/2}^{\pi/2} C_{\l}(\t) \t^2 = 
\frac{a(\l)}{2\pi}
\int_{-\pi/2}^{\pi/2} \t^2 (1+\cos (\t)) \cos^{2 \l}(\t) \, d\t   +
    \frac{a(\l) - b(\l)}{4\pi} \int_{-\pi/2}^{\pi/2} \t^2 \cos^{2 \l+1 }(\t) \, d\t }\\ &
= & \sqrt{\frac{\l}{\pi}} \int_{-\infty}^{\infty} \t^2 e^{-\l \t^2} \, d \t + O
\Bigl( \int_{-\infty}^{\infty} ( \l^{-1/2} \t^2   + \l^{1/2} \t^4 
    + \l^{3/2} \t^6) e^{-\l \t^2}  \,  d \t \Bigr)   
=  \frac{1}{2 \l}  + O \Bigl(\frac{1}{\l^{5/2}} \Bigr). 
\eeqn   
There remains the contribution due the last term on the right of
(\ref{equ:taylor}). The kernel $C_{\l}$ is to be integrated  
against
\be
\label{equ:nu}
\Bigl| \int_0^{\t} \Bigl( h^{\prime \prime}(\a)  -  h^{\prime
  \prime}(0) \Bigr) (\t - \a)  d \a \Bigr| \le \t^2 \hspace{-.1cm}
  \max_{-1 \le \a
  \le 1} |h^{\prime \prime}(\t \a) - h^{\prime \prime}(0)| :=
  \t^2 \tilde{\nu}(\t,h)
\ee
from $-\pi/2 \le \t \le \pi/2$.
By continuity, the nonnegative function $ \tilde{\nu}$ satisfies
$\tilde{\nu} \downarrow 0$ with  $\t \downarrow 0$.  Running through 
the above arguments will explain why the integral in question is bounded
by a constant multiple of 
\[
  \sqrt{\l} \int_{-\l^{-1/4}}^{\l^{-1/4}} \t^2 \tilde{\nu(\t,h)} e^{-\l \t^2} \, d \t 
   \le \frac{2}{\l} \max_{|\t| \le \l^{-1/4}} \tilde{\nu}(\t,h) := \frac{1}{\l} \nu(\l,h)
\]
with $\lim_{\l \uparrow \infty} \nu(\l,h) = 0$ as promised in (\ref{equ:er1}).
The comment following that display stems from the fact that if $h^{\prime \prime}$ is
H\"older continuous with exponent $\d > 0$, the right hand side of (\ref{equ:nu})
may be replaced by a constant times $\t^{2+\d}$.  The ensuing integral in the
last display would then decay like $\l^{-(1+\d/2)}$.

The details behind (\ref{equ:er2}) are much the same.
One may assume that $h(\t) = h_{+} \t^{+} - h_{-} \t^{-}$ in
some neighborhood of the origin and the leading 
order comes out of $\int_0^{\infty} \t e^{-\l \t^2} d \t =
1/2 \l$. The proof is finished.

\bigskip

\noindent{\bf Proof of Theorem \ref{thm:p1}}  As $f$ and $\tilde{g}$
each have one $L^2$ derivative,  their convolution $\phi(\t) = (f \ast 
\tilde{g})(\t)$ is twice continuously  differentiable. Recalling
the covariance formula
\[
   {\CV}_N \Bigl( X(f), X(g) \Bigr) = N \phi(0) - S_N(\phi)  - S_N^{\prime}(\phi)
\]
for 
\[
S_N (\phi)   =  \sum_{\l = 0}^{N-1}
                \int_{\T} C_{\l}(\t) \phi(\t) d \t 
\ \mbox{ and } \ 
S_N^{\prime}(\phi) = \sum_{N/2 \le \l \le {N}} C_{\l} \ast ( \phi - D_{2N -2\l -2} \phi) (0),
\]
we are in the setting of the first half of Lemma \ref{lem:error}. The 
statement (\ref{equ:er1}) implies that
\be
\label{equ:SN}
S_N (\phi)   =   N \phi(0) +
                 \frac{1}{4} \phi^{\prime \prime}(0) \sum_{\l = 1}^{N} \frac{1}{\l}
                  + \sum_{\l = 1}^N \frac{\nu(\l, \phi)}{\l} + O(1).
\ee
(Here, the $O(1)$ term includes both the sum of the order $\l^{-3/2}$ errors 
in (\ref{equ:er1}) as well as the sum of first several terms $(C_{\l} \ast \phi)(0)$
where the estimate is of no affect.)  The second term produces the $\log N$ 
figuring into the leading order growth of covariance.
Generally $\nu(\l, \phi) = o(1)$ and we can conclude only that
the third term $\sum_{\l = 1}^N \frac{\nu(\l, \phi)}{\l} = o(\log N)$.  The
point of the closing remark in the theorem's statement  is that
the proposed additional decay on $|\widehat{f}(k)|$ (or $|\widehat{g}(k)|$)
for $k \uparrow \infty$
implies H{\"o}lder continuity of  $\phi^{\prime \prime}$.
In particular,
for positive $\delta$ satisfying $\delta/2 < 1$,
\[
\Bigl|
\phi^{\prime \prime}(\t) -  \phi^{\prime \prime}(0) \Bigr|  \le 
\sum_{k = -\infty}^{\infty} k^2 | \widehat{f}(k)| | \widehat{g}(k)|
\Bigl| 1 - e^{i k \t} \Bigr|  \le  4 | \t |^{\d/2}
\sum_{k = -\infty}^{\infty} |k|^{2+\d/2}  | \widehat{f}(k)| | \widehat{g}(k)| 
\]
in which the last sum is finite by Schwartz's 
inequality.
By the related remark in the statement of Lemma \ref{lem:error} it follows
that the term  $\l^{-1} \nu(\l, \phi)$ is summable in this case.  In either
case, $N \phi(0) - S_N(\phi)$ already accounts for the advertised asymptotic
behavior of the covariance (note that $\phi^{\prime \prime}(0) = - \int_{\T} f^{\prime}(\t)
g^{\prime}(\t) d \t$).

The remaining (error) term
$S_N^{\prime}(\phi)$ is seen to be of 
constant order under the present smoothness assumptions. We
write
\be
\label{equ:Spr}
 S_N^{\prime}(\phi) 
                    = \sum_{N/2 \le \l \le {N}} \ \sum_{|k| > 2N - 2 \l -2} {\widehat{C}}_{\l}(k) \widehat{\phi}(k),
\ee 
and introduce the simple facts:
${\widehat{C}}_{\l}(k) = 0$ for $k > 2\l +1$, and, in general, 
 $| \widehat{C}_{\l}(k)| \le 2$. It follows that
\beqn
  | S_N^{\prime}(\phi)|  & \le & c_1 \sum_{0 \le \l \le N/2} \ \sum_{N - 2\l \le k \le N + 2 \l } 
                                   \Bigl(|\widehat{\phi}(k)|
          + |\widehat{\phi}(-k)| \Bigr)  \\ 
 & \le & c_2 \sum_{0 \le k \le N}  |k| \Bigl(|\widehat{\phi}(k)|
          + |\widehat{\phi}(-k)| \Bigr)  \le c_3 \Bigl( \sum_{k  = -\infty}^{\infty} |k|^4 |\widehat{\phi}(k)|^2 
             \Bigr)^{1/2} < \infty
\eeqn 
after changing variables and then the order of summation in inequalities one and two, an application of Cauchy-Schwartz 
for the third, and at last the given fact: $ \sum |k|^4 |\widehat{\phi}(k)|^2$  $\le$ $ (\sum |m|^2 |\widehat{f}(m)|^2)$
$(\sum |n|^2 |\widehat{g}(-n)|^2)$ $< \infty$.
The proof is finished.

\bigskip

\noindent{\bf Proof of Theorem \ref{thm:num2}}  When $f(\t) = g(\t) =
{\mathbbm 1}_{[-\a/2,\a/2]}(\t)$, the convolution $: =  \phi_{\a}(\t)$ 
is the tent function
\[
   \phi_{\a}(\t)  = \frac{1}{2 \pi} \int_{-\frac{\a}{2}}^{\frac{\a}{2}} 
                     {\mathbbm{1}}_{[-\frac{\a}{2}, \frac{\a}{2} ] }(\t^{\prime} - \t) \, d \t^{\prime} = 
\frac{1}{2\pi} \Bigl[ \a - |\t | \Bigr]
    {\mathbbm{1}}_{[-\a, \a]}(\t),
\]    
and the main contribution to the variance is then 
\be
\label{equ:phvar}
 N \phi_{\alpha}(0) - S_N(\phi_{\a})  =   N \frac{\a}{2 \pi} - \frac{\a}{4 \pi^2}
   \sum_{\l=0}^{N-1} \int_{-\a}^{\a} C_{\l}(\t) \, d \t
   +  \frac{1}{4 \pi^2} \sum_{\l =0}^{N-1} \int_{- \a}^{ \a} |\t| C_{\l} (\t)
       \, d \t. 
\ee  
Let us first consider this object for the
various situations, $\a$ fixed or $ \a \downarrow 0$ as $N \uparrow \infty$, 
and return to the error term $S_N^{\prime}(\phi_{\a}) $ at the end.

From the proof of Lemma \ref{lem:error} it follows
that for $\a > 0$ fixed $ \frac{1}{2\pi} \int_{- \a}^{ \a}  C_{\l}(\t) d\t$ $
= 1 + O(\l^{-3/2})$, and  so the first two terms in (\ref{equ:phvar}) contribute
something of constant order. On the other hand, from the conclusion of the second lemma we have
that
\[
     \frac{1}{4 \pi^2} \sum_{\l =0}^{N-1} \int_{- \a}^{ \a} |\t| C_{\l} (\t) d
     \t  = \frac{1}{2 \pi} \sum_{\l=1}^{N-1} \Bigl[ \frac{1}{\sqrt{\pi \l}}  +
     O \Bigl( \frac{1}{{\l}^{3/2}} \Bigr) \Bigr] + O(1) = \frac{1}{
     {\pi}^{3/2} } \sqrt{N} + O(1),
\]
which is the advertised result (\ref{equ:v3}) for order one regions of the phase.

If $\a \downarrow 0$ with $N \uparrow \infty$, we start in the
borderline case in which $ \a(N) = \b/\sqrt{N}$  with a fixed $\b > 0$.
At this scale the quantities in (\ref{equ:phvar}) go over into
approximating Riemann sums. First there is,
\beqn
\lefteqn{ \frac{1}{4 \pi^2}  \sum_{\l = 1}^{N}  \int_{-\a(N)
}^{\a(N)}  |\t| C_{\l}(\t) d \t  =  
\frac{\sqrt{N}}{2 \pi^{3/2}}  \times \frac{1}{N} \sum_{\l = 1}^{N} \sqrt{\frac{N}{\l}}  \int_{- \b
\sqrt{\frac{\l}{N}}}^{\b \sqrt{\frac{\l}{N}}}  |\t| e^{-\t^2} d \t  +
O(1)}  \\    & = &  \sqrt{N} \frac{1}{ \pi^{3/2}} \int_0^1 \frac{1}{\sqrt{x}} 
                    \int_{0}^{\b \sqrt{x}} \t 
                    e^{-\t^2} \, d \t d x  
 + O \Bigl( \sqrt{N} \sum_{\l = 1}^{N} \int_{\frac{\l-1}{N}}^{\frac{\l}{N}}
             \Bigl[ F_{\b}({\l}/{N} ) - F_{\b} ( x) \Bigr]  dx \Bigr) + O(1)
\eeqn 
for $F_{\b}(x) = \frac{1- e^{-\b x}}{\sqrt{x}}$, and that error term is controlled 
(for whatever $\b$) as in
\[
\sqrt{N} \sum_{\l = 1}^{N} \int_{\frac{\l-1}{N}}^{\frac{\l}{N}}
                \Big[ F_{\b}({\l}/{N} ) 
                      - F_{\b} ({x} ) \Bigr]  dx  \\
 \le   \frac{1}{\sqrt{N}} \int_0^1 |F_{\b}^{\prime} (x)| dx = O\Bigl( \frac{1}{\sqrt{N}} \Bigr).
\] 
(The difference between an integral and its Riemann sum being less than $1/N \times$  
the total variation of the integrand.)
In a similar way we find that
\[
   N \frac{\a(N)}{2 \pi} - \frac{\a(N)}{4 \pi^2}
   \sum_{\l=0}^{N-1} \int_{-\a(N)}^{\a(N)} C_{\l}(\t) \, d \t
= \sqrt{N} \Bigl[ \frac{\b}{2 \pi} - \frac{\b}{2 \pi} \int_0^1 
              \int_{-\b \sqrt{x}}^{\b \sqrt{x}} \frac{e^{-\t^2}}{\sqrt{\pi}} \, d \t dx \Bigr] + O(1). 
\]
Combined, the last three displays translate as,
\beq
\label{equ:Iarg}
\lefteqn{ N \phi_{\frac{\b}{\sqrt{N}}} (0) - S_N( \phi_{\frac{\b}{\sqrt{N}}}) }  \\
  & = & \sqrt{N} \frac{1}{\pi^{3/2}} 
        \Bigl[  \int_0^1 \frac{1- e^{-\b x^2}}{2 \sqrt{x}} \, dx 
         + \b \int_0^1 \int_{\b \sqrt{x}}^{\infty} e^{-\t^2} \, d \t d x                                      
        \Bigr] + O(1)  \nonumber \\
  & := & \sqrt{N} \frac{1}{\pi^{3/2}} I_{arg}(\b) + O(1). \nonumber
\eeq

Next, if $\sqrt{N} \a(N) \uparrow \infty$ one can clearly apply the same steps
with the point of view that $\beta = \beta(N)
\uparrow \infty$.  Noting that 
\[
  I_{arg}(\beta) = 1 - \frac{c_1}{\b^{1/4}} + \frac{c_2}{\beta}
\]
as $\b \uparrow \infty$ completes the explanation of that case.
Finally, if $\sqrt{N} \a(N) \downarrow 0$ while maintaining $N \a(N) \uparrow \infty$,
things are a bit different. In the range $\l \le N$ and $\sqrt{ N} \t = o(1)$,
\[
   C_{\l}(\t) = 2 \sqrt{ \pi \l} \Bigl( 1 - \l \t^2 + O(\l^2 \t^4) 
                  \Bigr) \Bigl(1 + O \Bigl(\frac{1}{\l} \Bigr) \Bigr)  
                    + O \Bigl( \frac{1}{\l^{3/2}} \Bigr)
\]
Substituting the above into (\ref{equ:phvar}) we
find that
\[
  \frac{\a}{4 \pi^2} \sum_{\l=0}^{N-1} \int_{-\a}^{\a} C_{\l}(\t) d \t 
 - \frac{1}{4 \pi^2} \sum_{\l=0}^{N-1} \int_{-\a}^{\a} |\t| C_{\l}(\t) d \t
  = 
   N \a \Bigl[  \frac{ \sqrt{N} \a}{3 \pi^{3/2}} - \frac{ (\sqrt{N} \a)^3}{ 30 \pi^{3/2}}
         \Bigr] + O(1),
\]
which anyway is $o(N \a)$,
identifying the leading order in this case as $N (\a /2  \pi) $.

To complete the proof,  we need to track the growth of $S_N^{\prime}$.
A bit of good fortune is that one may compute ${\widehat{\phi}}_{\a}(k)$
$= \frac{\sin^2( k \a)}{\pi k^2} \ge 0$ which leaves us to control
\be
\label{equ:Spr2}
 |S_N^{\prime}(\phi_{\a})| \le c_1 \sum_{N/2 \le \l \le N} \ \sum_{ 2N - 2 \l \le |k| \le 2 \l} 
|\widehat{C}_{\l}(k)| \, \Bigl(\frac{\sin^2(k \a)}{k^2} \Bigr),
\ee
see (\ref{equ:Spr}). Simply bounding $|\widehat{C}_{\l}(k)| \sin^2(k \a)$ by a constant
will produce the estimate $S_N^{\prime}(\phi_{\a}) = O(\log N)$ stated for the cases
where $\liminf_{N \uparrow \infty} \sqrt{N} \a > 0$.

If however, $\a = o(1/\sqrt{N})$ we need to, and can, do better.
Bring in the additional estimate $|\widehat{C}_{\l}(k)| \le c_2 \sqrt{\l} |k|^{-1}$
courtesy Lemma {\ref{lem:error}} which has bite for $|k| >> \l$.
With this in mind the sum is split as follows: with
 $\l^*(N) \simeq N - \frac{1}{2} \sqrt{N}$
the point where $2N -2\l = \sqrt{\l}$ over $\l > N/2$,
\[
| S_N^{\prime}(\phi_{\a})| \le \Bigl[ 
  \sum_{{ N/2 \le \l \le \l^*(N)} \atop 
                   {2N - 2 \l \le |k| \le 2 \l}}  +   
\sum_{{ \l^*(N) \le \l \le N } \atop {
                   \sqrt{\l} \le |k| \le 2 \l}}   +
\sum_{{\l^*(N) \le \l \le N } \atop {  
                   1 \le |k| \le \sqrt{\l}}}
   \Bigr]  \,  |\widehat{C}_{\l}(k)| \, \Bigl(\frac{\sin^2(k \a)}{k^2} \Bigr)
   :=  {\cal A}_N + {\cal B}_N + {\cal C}_N.  
\]
Moving from left to right, and first bounding the summand by
$\sqrt{\l} |k|^{-3}$ we have  
 \[
  {\cal A}_N \le c_3  \sum_{N/2 \le \l \le N - \frac{1}{4} \sqrt{N}} 
   \frac{\sqrt{\l}}{ (N -  \l)^2 } \le {c_4 } {\sqrt N} \hspace{-.1cm}
          \sum_{\frac{1}{4} \sqrt{N} \le \l \le  N/2}  \frac{1}{\l^2} = O(1), 
\]
as well as 
\[
{\cal B}_N \le c_5 
  \sum_{{N - \sqrt N} \le \l \le N} \, \sum_{\sqrt{\l} \le k \le \infty} 
   \frac{\sqrt{\l}}{k^3}
     \le c_6 \sum_{N- {\sqrt N} \le \l \le N}  \frac{1}{\sqrt \l} = O(1).
\]
For  ${\cal C}_N$ we finally use the fact that $\a \downarrow 0$: 
as 
$k^{-2} \sin^2(\a k) \le 2 \a^2$  for $\a = o(1/{\sqrt N})$ and $k \le {\sqrt N}$
we conclude
\[
 {\cal C}_N \le c_7 \, {\a}^2  \hspace{-.2cm} \sum_{N - \sqrt{N} \le \l \le N } \sqrt{\l}
    \le c_8 \, (N \a^2) = o(1) 
\]
as needed.  In closing 
we note that a careful review of the ${\cal C}_N$ will demonstrate
that the stated $O(\log N)$ error term in the cases  where
$\liminf_{N \uparrow \infty} \sqrt{N} \a 
> 0$ cannot be improved.  The proof is finished.

\bigskip

\noindent{\bf Proof of Theorem \ref{thm:lcov} (angular case)}
In each on the four cases, the covariance of the pair  
$(\# \Theta {[\a, \b]} , \# \Theta [\d,\gamma] )$  
is connected to a $\phi_{k}(\t)$ ($k = 1,2,3,4$) 
which is either constant in a neighborhood
of the origin or possesses a corner at that point.  
In particular, when the intervals in question are either disjoint or else $[\a, \b]$ sits 
inside of $[\d, \gamma]$ we have 
\[
   \phi_{1,2}(\t) := \frac{1}{2\pi}
       \int_{\d - \t}^{\gamma - \t} {\mathbbm 1}_{[\a, \b]}(\t^{\prime}) \, d \t^{\prime} 
  = \left\{
    \begin{array}{ll} 0 & \mbox{ on } [ \b -\d, 2\pi + \a - \gamma]  \mbox{ when  } \a < \b < \d < \gamma \\
                      \frac{\b - \a}{2 \pi} 
                      & \mbox{ on } [ \d - \a, \gamma - \b ] \mbox{ when } 
                      \d < \a < \b < \gamma
    \end{array} \right.                       
\]
The relevant computation is then: for 
$\phi(\t) = \phi(0)$ over  $ - \ep \le \t \le \ep$, 
\[
  N \phi(0) - \sum_{\l =0}^{N-1} \int_{\T} C_{\l}(\t) \phi(\t) \, d\t
   = N \phi(0) - \phi(0) \sum_{\l=0}^{N-1} \Bigl[ 
 \int_{- \l \ep }^{\l \ep} e^{-\t^2} \frac{d \t}{\sqrt{ \pi}} + O(\l^{-3/2}) \Bigr]
   = O(1).
\]
As for $S_N^{\prime}(\phi_{1,2})$,
we will estimate  $f_m^{1,2}(\t) := (\phi_{1,2} - D_m \ast \phi_{1,2})(\t)$ (which is evaluated
at $m = 2N - 2l -2$ and then summed over $\l$) near the origin. In particular,
we have already seen that by adding a constant if necessary we can assume that 
 $\phi_{1,2}(\t) = 0$ in $[-\ep,\ep]$, and we want to show that $f_{m}$ can be made
uniformly small throughout say $[-\ep/2, \ep/2]$. Now,
\[
 f_m^{1,2}(\t) = \int_{T} \frac{\phi_{1,2}(\t - \t^{\prime})}{\sin( \frac{1}{2} \t^{\prime})} 
                 \sin \Bigl( (m + \frac{1}{2} \Bigr) \t^{\prime}) \, d \t^{\prime} : =  \int_{T} g_{1,2}(\t^{\prime}, \t)  
                     \sin \Bigl( (m + \frac{1}{2}) \t^{\prime} \Bigr) \, d \t^{\prime}, 
\]  
and if $|\t| \le \ep/2$ then $\phi_{1,2}(\t - \t^{\prime}) = 0$ for all $|\t^{\prime}| \le \ep/2$.
As otherwise $\phi_{1,2}$ is piece-wise linear, the second derivative of $g_{1,2}(\t^{\prime})$ will certainly 
be integrable away from  $[-\ep/2, \ep/2]$.  Therefore, integrating by parts twice 
will show that $f_{m}^{1,2}(\t) = O(m^{-2})$ throughout $|\t| \le \ep/2$. This in turn  will imply that
\[
 \Bigl( C_{\l} \ast f_{2 N - 2 \l - 2}^{1,2} \Bigr)(0) = O \Bigl( \frac{1}{(N -\l)^2} \Bigr) + O( \l^{-3/2} )  
\]
which is summable over the appropriate range of $\l$:  $S_N^{\prime}(\phi_{1,2})$  is 
of constant order.  The de-correlation at the proposed rate is thus
established.

Taking next the situation that $\a < \b = \d < \gamma$,  it is enough 
to note that there exists an $\ep > 0$ with
\[
  \phi_{3}(\t) := \frac{1}{2 \pi}
 \int_{\a - \t}^{\b - \t} {\mathbbm 1}_{[\b, \gamma]}(\t^{\prime}) \, d \t^{\prime}
     =   \left\{
    \begin{array}{ll} \frac{-\t}{2 \pi}  & \mbox{ for }  \t \in [-\ep, 0] \\
                       0   & \mbox{ for }  \t \in [0, \ep], 
    \end{array} \right.  
\] 
whereas if $\a = \d$ and $\b \neq \gamma$ (take $ \b < \gamma$ for convenience) 
we may say that
\[
\phi_{4}(\t) := \frac{1}{2 \pi} 
   \int_{\a - \t}^{\b - \t} {\mathbbm 1}_{[\a, \gamma]}(\t^{\prime}) \, d \t^{\prime}
     =   \left\{
    \begin{array}{ll}  \frac{\b - \a}{2 \pi}   & \mbox{ for }  \t \in [-\ep^{\prime}, 0] \\
                       \frac{\b - \a - \t}{ 2 \pi}  
                  & \mbox{ for }  \t \in [0, \ep^{\prime}] 
    \end{array} \right.  
\]
for some $\ep^{\prime} > 0$.
It will then follow from the second half of  Lemma {\ref{lem:error}} that
\[
   - \sum_{\l=0}^{N} \int_{\T} C_{\l}(\t) \phi_3(\t) \, d \t =  + \frac{1}{2} \sqrt{ \frac{N}{\pi^3}} 
 + O(1) \ \ \mbox{ and } \ \ 
   -   \sum_{\l=0}^{N} \int_{\T} C_{\l}(\t) \phi_4(\t) \, d \t =  - \frac{1}{2} \sqrt{ \frac{N}{\pi^3}} 
 + O(1).
\]
Furthermore, re-tooling that part of the previous proof (for the variance estimates) 
related to the $S_N^{\prime}$ term will show that $S_N^{\prime}(\phi_{3,4}) = O(\log N)$;
the corner at the origin being the cause of this increased growth.
The proof is finished.

\section{Radial Statistics}
\setcounter{num}{4}
\setcounter{equation}{0}
\label{sec:rad}

We begin with the proof of Theorem {\ref{thm:r1} which relies on
the product formula
\be
\label{equ:meanfactor}
   E \Bigl[ \prod_{k=1}^N f(|z_k|) \Bigr]  = \prod_{k=1}^N E \Bigl[  f \Bigl( \sqrt{ \frac{1}{N} s_{k} } \Bigr) \Bigr].
\ee
where you will recall that $s_k = \eta_1 + \cdots + \eta_k$ with the $\eta$'s independent exponentials of mean one.
While related to the covariance formula (\ref{equ:cov1}), it is easier to see (\ref{equ:meanfactor}) as direct
consequence of (\ref{equ:comrad}) derived in the next section.

\bigskip

\noindent{\bf Proof of Theorem \ref{thm:r1}}  Let $h$ satisfy the assumed boundedness and regularity. For convenience we 
consider test the function $r \rightarrow h(r^2)$ with a
Gaussian limit for $\sum_{k=1}^N h(|z_k|^2) - E [ \sum_{k=1}^N h(|z_k|^2) ]$ proved by 
computing Laplace transforms. 
In particular, the formula (\ref{equ:meanfactor}) is applied with 
$f(r) = \exp{[ \lambda h(r^2) ]}$ and $\lambda$ a parameter lying in a fixed neighborhood of the
origin to find that
\beq
\label{equ:fir}
  E \Bigl[ \exp{ \Bigl\{ \lambda  [ X(h) - E X(h)  ]  \Bigr\} } \Bigr] & = & 
  E \Bigl[  \exp{ \Bigl\{ \lambda [ \sum_{k=1}^N ( h(|z_k|^2)  - E h(|z_k|^2) ) ] \Bigr\} } \Bigr] \\
& = & \prod_{k=1}^N  
E \Bigl[ \exp{ \Bigl\{ \lambda[  h ({s_k}/{N}) - E h ({s_k}/{N} ) ] \Bigr\} } \Bigr]. \nonumber
\eeq 
We wish to expand each appearance of $h(s_k/N)$ as 
in $h(s_k/N) = h(k/N)  + h^{\prime}(k/N) ((s_k - k)/N) + \frac{1}{2} h^{\prime \prime}(k/N)
( {(s_k - k)}/N )^2 + R_h(s_k/N, k/N)$.  Toward this, make the definitions 
$\bar{h}(s_k/N) = h(s_k/N) - E[h(s_k/N)]$,
\[
  {\cal M}_{N,k}(h) := \exp{  \Bigl\{ \lambda h^{\prime}({k}/{N} ) \Bigl( \frac{s_k - k}{N} \Bigr) \Bigr\}},
\] 
and 
\[
 {\cal R}_{N,k}(h) := \exp{  \Bigl\{ \ \frac{\lambda}{2N^2}  h^{\prime \prime }({k}/{N} ) 
         \Bigl( ( {s_k - k} )^2 - k \Bigr)
                                       + \lambda \bar{R}_h( s_k/N, k/N )  \Bigr\}}
\]
in which $\bar{R}_h( s_k/N, k/N) = 
R_{h}(s_k/N, k/N) - E [R_h(s_k/N, k/N)]$, and the right hand side
of (\ref{equ:fir}) may be continued as the product of
\beq
\label{equ:sec}
\lefteqn{E \Bigl[ \exp{ [ \lambda \bar{h}(s_k/N) ] } \Bigr] = 
 E \Bigl[ {\cal M}_{N,k}(h) \Bigr] } \\
&  &  +  E \Bigl[ \Bigl( e^{\lambda  \bar{h}(s_k/N)} -  {\cal M}_{N,k}(h) \Bigr) 
                              {\mathbbm 1}_{ \{\frac{s_k}{N} \ge 1+ \ep \} }      
                          \Bigr] 
 + E \Bigl[    {\cal M}_{N,k}(h) \Bigl(  {\cal R}_{N,k}(h) - 1 \Bigr) 
                           {\mathbbm 1}_{\{ \frac{s_k}{N} \le 1+ \ep \} }  \Bigr] \nonumber
\eeq
from $k = 1$ to $N$.

Computing the product of the first term in the right of (\ref{equ:sec}) produces the Gaussian structure: 
choosing $\lambda$ so that $|\lambda| \max_{r \in [0,1]} |h^{\prime}(r)| < 1/2$,
\beqn 
\lefteqn{ \log \prod_{k=1}^N E \Bigl[  {\cal M}_{N,k}(h) \Bigr]
          = - \sum_{k=1}^N \Bigl(  \lambda h^{\prime}({k}/{N} ) \frac{k}{N} 
                   + k \log \Bigl( 1 - \frac{\lambda h^{\prime}({k}/{N} )}{N} \Bigr)
              \Bigr) } \\ 
& = &  - \sum_{k=1}^N k \Bigl( \lambda h^{\prime}({k}/{N} ) \frac{1}{N} \Bigr)^2 + O \Bigl( \sum_{k=1}^N \frac{k}{N^3} \Bigr)
            = - \frac{1}{2} \lambda^2 \int_0^1 ( h^{\prime}(r))^2 \, r dr  + o(1).
\eeqn 
Now similar 
estimates to those used in the last display will show that $E \Bigl[  {\cal M}_{N,k}(h) \Bigr]$
is bounded above and below by positive constants independent of $k$ or $N$, and so the proof is completed
by establishing that
\[
 \sum_{k=1}^N  \Bigl\{ E \Bigl[ \Bigl( e^{\lambda  \bar{h}(s_k/N) } +  {\cal M}_{N,k}(h) \Bigr) 
                              {\mathbbm 1}_{\{\frac{s_k}{N} \ge 1+ \ep \}}
                          \Bigr] 
 + E \Bigl[  {\cal M}_{N,k}(h) (  {\cal R}_{N,k}(h) - 1) {\mathbbm 1}_{ \{ \frac{s_k}{N} \le 1+ \ep \} }  \Bigr]  
               \Bigr\} \rightarrow 0
\] 
as $N \uparrow \infty$.

For the first part, let $c_1 = \sup_{0 \le r < \infty} h(r)$ and write
\beqn
\lefteqn{ \hspace{-2.5cm}
E \Bigl[ \Bigl( e^{\lambda  \bar{h}(s_k/N)} +  {\cal M}_{N,k}(h) \Bigr) 
                              {\mathbbm 1}_{\{\frac{s_k}{N} \ge 1+ \ep \}}
                          \Bigr] 
 \le  e^{2 |\lambda| c_1} P \Bigl( \frac{s_k}{N} \ge 1 + \ep \Bigr) 
         +   E \Bigl[ e^{s_k/2N}  {\mathbbm 1}_{\{\frac{s_k}{N} \ge 1+ \ep\}} \Bigr] } \\
  & \le & e^{2 |\lambda| c_1} P \Bigl( \frac{s_k}{N} \ge 1 + \ep \Bigr) +  \Bigl( \frac{N}{N-1} \Bigr)^{k/2} 
          \sqrt{  P \Bigl( \frac{s_k}{N} \ge 1 + \ep \Bigr)} 
\eeqn
in which the last line is easily bounded by a constant multiple of $ e^{-\gamma N} $
for  $\gamma > 0$ depending on $\ep$ but not on $k$.  For the remaining part, note first that
on $\{ s_k/N \le 1 + \ep \}$ we have that ${\cal M}_{N,k}(h) \le \exp{ (1 + \ep)}$ and also that
the exponent in ${\cal R}_{N,k}(h)$ may be bounded uniformly in $k$ and $N$. Since 
$| e^{a} - 1| \le |a| e^{|a|}$ we have that
\beq
\label{equ:lassplit}
\lefteqn{ \hspace{-2cm} 
E \Bigl[  {\cal M}_{N,k}(h) | {\cal R}_{N,k}(h) - 1 | {\mathbbm 1}_{ \{ \frac{s_k}{N} \le 1+ \ep \} }  
 \Bigr]} \\  
&  &  \le c_2    \frac{k}{N^2} \, E \Bigl[ | \Bigl( \frac{s_k - k}{\sqrt{k}} \Bigr)^2 - 1 | \Bigr] 
           +  c_3 E \Bigl[|  {R}_h( s_k/N , k/N) | \Bigr] \nonumber 
\eeq  
with constants $c_2$ and $c_3$.  
Clearly $E[ | {(\frac{s_k -k}{\sqrt{k}})}^2 - 1 |] = o(1)$  for $k \uparrow \infty$.  It follows
that the first term on the right of (\ref{equ:lassplit}), summed from $k = 1$ to $N$, is also $o(1)$
as $N \uparrow \infty$.  For the final term in (\ref{equ:lassplit}),
by the regularity of $h$ we may assume that
$ |\bar{R}_h(s_k/N, k/N)| \le c_4 |\frac{s_k - k}{N}|^{2 +\alpha}$ (having applied this estimate while 
the restriction $\{ s_k/N \le 1 + \ep\}$ was still in place). The upshot is that 
\[
\label{equ:rem2}
 \sum_{k=1}^N E \Bigl[ | R_h( s_k/N, k/N)|  \Bigr] \le  \frac{c_4}{N^{2 + \alpha}}  
  \sum_{k=1}^N E \Bigl[ |s_k - k |^{2 +\alpha} \Bigr]  
   \le c_5  \frac{1}{N^{2 + \alpha}} \sum_{k=1}^N k^{\frac{2 + \alpha}{2}} \simeq N^{-\a/2} 
\]
as desired. Here, the second inequality may be arrived at rather directly in the present setting of
exponential random variables. More generally, an   
inequality due to Marcinkiewicz and Zygmund (which may be found in \cite{CT}) will explain 
why $ E [ | \sum_{\l =1 }^k  x_{\l} |^p ] \simeq k^{p/2}$  for $p > 2$  and $\{ x_{\l} \} $ independent copies of a mean-zero $x_1$  
for which $E[|x_1|^p] < \infty$.

Finally, the statement concerning the asymptotic behavior of the covariance may be gleaned from the 
above by setting $h = \beta f + \gamma g$
for parameters $\b$ and $\gamma$.  Having controlled moment generating functions, we have adequate tightness
to draw conclusions on the limiting covariance as well as higher moments.  The proof is finished.
\bigskip

Moving on to the covariance of the number statistic $\# {\I}{[ \cdot,\cdot ]}$,
the basic error estimate required has in fact already 
been done for us. It is embodied in the classical Edgeworth
expansion providing corrections to the local central limit
theorem.  The statement is: 
with $p(a,M)$ the density of the random variable 
$\frac{1}{\sqrt{M}} \sum_{\ell = 1}^M (\eta_{\ell} - 1)$ at $a$,
\begin{equation}
\label{equ:edgeworth}
\sup_{-\infty < a < \infty} \Bigl| p(a,M) - \frac{1}{\sqrt{2 \pi}} e^{-a^2/2}
                        - \frac{1}{\sqrt{M}}  \frac{1}{\sqrt{\pi/2}}                a^3 e^{-a^2/2}
                            \Bigr| = O(M^{-1}),
\end{equation}
see for example  \cite{BR}, Corollary 19.4.

\bigskip

\noindent{\bf Proof of Theorem \ref{thm:num1}}  
Start with a fixed interval $[a,b]$ ($0 < a \le b < 1$)
of the moduli for $N \uparrow \infty$.  Denote by  $p_{N,k}(A) 
= P(\frac{1}{N} s_k \in A)$.
The variance
of $\# \I {[a,b]} $ then reads
\be
{\V}_N \Bigl(  \# {\I}{[a,b]} \Bigr) 
= \sum_{k=1}^N p_{N,k}([a^2,b^2]) \Bigl( 1 - p_{N,k}([a^2, b^2]) \Bigr)
:=  \sum_{k=1}^{N} p_{N,k}([a^2,b^2]) p_{N,k}([a^2, b^2]^C),  
\ee
and the (classical) central limit theorem explains why it should be 
only $O(\sqrt{N})$ neighborhoods of $k = N a^2$ and
$k =N b^2$ that contribute to the sum on the right.  
Breaking up the sum in accordance with that observation,
we write:
\[
 {\V}_N \Bigl( \# {\I} {[a,b]} \Bigr) 
    = \left[ \sum_{ |k- Na^2| \le \d_N \sqrt{N}}
    + \sum_{|k - Nb^2| \le \d_N \sqrt{N}}  
    + \ {\sum}^{\prime} \  \right] \, p_{N,k}([a^2,b^2]) 
      p_{N,k} ([a^2,b^2]^C). 
\]
Here ${\sum}^{\prime}$ has in general three components ($k \in$
$[0, Na^2 - \d_N \sqrt{N}]$ $\cup [Na^2 + \d_N \sqrt{N},
Nb^2 - \d_N \sqrt{N}] $ 
$ \cup [Nb^2 + \d_N \sqrt{N}, N]$), and $\delta_N$ is some
rate chosen so that this sum is negligible for $N \uparrow \infty$. 
We spell out the estimate for the first component, the rest 
being nearly identical. 

With $\{ \lambda_k \}$ a sequence of positive constants
less than one to be determined, we have:
\beq
\label{equ:tail}
\lefteqn{
\sum_{k \le Na^2 - \d_N \sqrt{N}} 
     p_{N,k}([a^2,b^2]) p_{N,k}([a^2,b^2]^C) 
  \le \sum_{k \le Na^2 - \d_N \sqrt{N}} P( s_k \ge N a^2) }  
   \\ 
& \le & \sum_{\d_N \sqrt{N} \le k \le N a^2}
      \exp{ \Bigl[ -N \Bigl( a^2 \lambda_k + (a^2 - \frac{k}{N}) 
                             \log(1-\lambda_k) \Bigr) \Bigr]} 
   \nonumber \\ 
& \le & \sum_{\delta_N \sqrt{N} \le k \le N a^2} e^{-k^2/2N a^2}
  \le c_1 \, \sqrt{N} 
  \int_{\delta_N}^{\infty} e^{-c^2/2a^2} \, dc 
  \le c_2 \frac{\sqrt{N} e^{-\delta_N^2/2a^2} }{ \delta_N}. 
   \nonumber
\eeq
The second inequality invokes Chernov's bound plus 
the change of variables $k \rightarrow \lfloor N a^2 \rfloor -k$. The third follows from the choice of
$\lambda_k = \frac{k}{Na^2}$, and the fourth from a Taylor expansion.  The
conclusion of the last line is that letting $\delta_N
= \sqrt{ \log N}$ implies $ \sum^{\prime} p_{N,k}(1 - p_{N,k})
$ $ = o(1)$.

Turning to the sums of
consequence (over $|k - Na^2| \le \d_N \sqrt{N}$ and 
$|k - Nb^2| \le \d_N \sqrt{N}$), it is evident that 
they will have equal contributions.  Consider the
first sum and rewrite the typical appearance of
$p_{N,k}([a^2,b^2])$ as: with $k = Na^2 +\l$ and
$\l$ ranging between $\pm \d_N \sqrt{N}$, 
\beq
\label{equ:term}
\lefteqn{ P \Bigl( \frac{1}{N} s_{Na^2 + \l} \in [a^2, b^2] \Bigr) }
 \\ & = & P \Bigl( \frac{1}{\sqrt{ \lfloor Na^2 + \l \rfloor }}
\sum_{1 \le m \le Na^2 +
 \l}  (\eta_m -1) \ge  -  \frac{\l}{\sqrt{ N}a}   + \ep_N(\l,a) \Bigr)
 - O \Bigl(e^{-c_{a,b} N} \Bigr).   \nonumber  
\eeq 
Here 
\[
  |\ep_N(\l,a)| \le 
\Bigl| \frac{\l}{\sqrt{Na^2}} - \frac{\l}{\sqrt{Na^2 + \l}} \Bigr| 
                + \frac{1}{\sqrt{N}} \Bigl( Na^2 + \l - \lfloor N a^2 + \l \rfloor \Bigr) 
            \le c_1 \frac{1}{\sqrt{N}} + c_2 \frac{\l^2}{N^{3/2}},
\]  
and the (exponentially small) error term stems
from a standard large deviation estimate. Applying
(\ref{equ:edgeworth}), we may then continue  (\ref{equ:term}) as in
\beq 
\label{equ:term1}
\lefteqn{p_{N,N a^2 +\l}([a^2,b^2])} \\
& = & \frac{1}{\sqrt{2 \pi}}
 \int_{-\frac{\l}{\sqrt{N} a} + \ep_N }^{\infty}  e^{-x^2/2} dx     
  + {\sqrt{\frac{2}{\pi  N}}}  \int_{-\frac{\l}{\sqrt{N} a }  +\ep_N }^{\infty} 
  x^3 e^{-x^2/2} dx   + O \Bigl(
 \frac{\log N}{N} \Bigr)  
 \nonumber \\ 
& = & \frac{1}{\sqrt{2 \pi}}
 \int_{-\frac{\l}{\sqrt{N} a }}^{\infty}  e^{-x^2/2} dx     
  + \sqrt{\frac{2}{\pi  N}}  \int_{-\frac{\l}{ \sqrt{N} a}}^{\infty} 
  x^3 e^{-x^2/2} dx   + O \Bigl( \frac{\l^2}{N^{3/2}} e^{-\l^2/N}  \Bigr)
  + O \Bigl(\frac{\log N}{N} \Bigr). \nonumber 
\eeq 
The  first equality is achieved by first cutting down the 
integration to  $|c| \le \log N$ and then
applying the Edgeworth expansion. The error produced by this process of 
limiting the range of integration and later restoring it up to $+\infty$ may
be controlled by the same large deviation procedure used already twice
before and the Gaussian estimate $\int_{a}^{\infty} e^{-x^2/2} dx \le
a^{-1} e^{-a^2/2}$.  The second equality is self-explanatory.
 
The leading order of the variance is then read off by subtracting the
square of (\ref{equ:term1}) and summing in $\l$:
\beq
\label{equ:rr}
\lefteqn{ \sum_{-\sqrt{N} \d_N \le \l \le \sqrt{N} \d_N}
    p_{N,N a^2 + \l }([a^2,b^2])  (1- p_{N,Na^2 + \l}([a^2,b^2]) )} \\
& = &  
    \sqrt{N} \hspace{-.2cm} \sum_{-\sqrt{N} \d_N \le \l \le \sqrt{N} \d_N}
      \frac{1}{\sqrt{N}} \Bigl\{ \int_{\frac{\l}{{\sqrt N} a}}^{\infty} 
             e^{-x^2/2} \frac{dx}{\sqrt{2 \pi}}  - 
            \Bigl[ \int_{\frac{\l}{{\sqrt N} a}}^{\infty} 
             e^{-x^2/2} \frac{dx}{\sqrt{2 \pi}} \Bigr]^2
     \Bigr\}  + {\cal E}_1(N) \nonumber \\
& = & \sqrt{N} \int_{-\infty}^{\infty} \Bigl\{ \int_{\frac{x^{\prime}}{a}}^{\infty} 
             e^{-x^2/2} \frac{dx}{\sqrt{2 \pi}}  - 
            \Bigl[ \int_{\frac{ x^{\prime}}{a}}^{\infty} 
             e^{-x^2/2} \frac{dx}{\sqrt{2 \pi}} \Bigr]^2 \Bigr\} dx ^{\prime} 
                   + {\cal E}_1(N) + {\cal E}_2(N). \nonumber
\eeq
The error involved in taking the sum over to the (full) integral
in the second line is easily controlled: 
\beq
\label{equ:samp}
{\cal E}_2(N) & \le & c_3 \sqrt{N} \int_{\d_N}^{\infty} \int_{x^{\prime}}^{\infty} e^{-x^2/2} dx dx^{\prime} 
                 + c_4 \sqrt{N} \hspace{-.15cm}  \sum_{|\l| \le \sqrt{N} \d_N}
                           \int_{\frac{\l}{{\sqrt{N}}}}^{\frac{\l+1}{{\sqrt N}}}
                 \Bigl| \int_{\frac{\l}{{\sqrt{N}}}}^{x^{\prime}} e^{-x^2/2} dx \Bigr| dx^{\prime} \\
              & \le & c_5 \frac{\sqrt{N}}{\d_N^2} e^{- \d_N^2/2} + c_6  \frac{1}{\sqrt{N}}
     \sum_{|\l| \le \sqrt{N} \d_N}  e^{- \l^2/2N} = o(1). \nonumber
\eeq
It remains to show that ${\cal E}_1(N) = O(1)$.  Terms three and four
in (\ref{equ:term1}) and their squares, summed over
the range $\pm \sqrt{N} \log N $ clearly remain finite as 
$N \uparrow \infty$. As for the second term, note that it
is already of the form $N^{-1/2}$ times an integrable 
function of $\l N^{-1/2}$.  That is, it is already properly
normalized so that the corresponding sums go over into
(convergent) integrals as $N \uparrow \infty$: an instance of this being,
\be
\label{equ:samp1}
   \frac{1}{\sqrt{N}} \sum_{|\l| \le  \sqrt{N} \d_N}  \int_{\frac{\l}{\sqrt{N}}}^{\infty} x^3 e^{-x^2/2} 
   \le  c_7 \frac{1}{\sqrt{N}}  \sum_{|\l| \le  \sqrt{N} \d_N}
 (1 + \frac{\l^2}{N} ) e^{- \l^2 / 2 N} \le c_8 \int_{-\infty}^{\infty} (1 + x^2) e^{-x^2/2} dx
\ee
after an integration by parts. 
Finally,
the integral in line three of (\ref{equ:rr}) is evaluated
as $ = \frac{a}{\sqrt{\pi}}$, and adding the corresponding 
term for $k \simeq Nb^2$ completes the proof of (\ref{equ:v1}).

The proof for mesoscopic scales (\ref{equ:v2}) builds on the above.
Let us set, $\sqrt{N} (b-a) = c$, and consider first the cases
that $c$ is fixed (the critical case) or $c \uparrow \infty$ with
$N$.  Also, we take the interval in question to be  $[a, a + c/\sqrt{N}]$
for a fixed $a$.  While one may also consider symmetric intervals shrinking
to $a$ at $N = \infty$, the proof is more in lines with the above if we
make the present convention.   In particular, retracing those arguments through
the first line of (\ref{equ:rr}),
it is evident that 
\[
  {\V}_N \Bigl(
  \# \I \Bigl[a , 
       a + \frac{c}{\sqrt{N}} \Bigr] \Bigr)
   =  \hspace{-.15cm} 
\sum_{|\l| \le \sqrt{N} \d_N}
      \Bigl\{ \int_{\frac{\l}{{\sqrt N} a}}^{2c + \frac{\l}{\sqrt{N} a}} 
             e^{-x^2/2} \frac{dx}{\sqrt{2 \pi}}  - 
            \Bigl[ \int_{\frac{\l}{{\sqrt N} a}}^{2c + \frac{\l}{\sqrt{N}a }} 
             e^{-x^2/2} \frac{dx}{\sqrt{2 \pi}} \Bigr]^2
     \Bigr\}  + {\cal E}_1(c,N).
\]
The error ${\cal E}_1(c,N)$ plays the same role as ${\cal E}_1(N)$ in (\ref{equ:rr})
and is amenable to the same arguments. While certain aspects such as the estimate
in (\ref{equ:samp1}) are now more elaborate due to the contribution of the upper
limit in the integral, everything goes through just as easily if one keeps in mind
that $c$ is at most $o(\sqrt{N})$. Next the outer sum is replaced by an integral
after an estimate nearly identical to (\ref{equ:samp}).  That is, the previous
display may be continued, and the variance given by
\beq
\label{equ:Imod}
         \frac{ \sqrt{N} a}{\sqrt{\pi}} \times 
  \sqrt{\pi}  \int_{-\infty}^{\infty}  \Bigl\{ \int_{x^{\prime}}^{2c +x^{\prime}} 
           \hspace{-.1cm}  e^{-x^2/2} \frac{dx}{\sqrt{2 \pi}} 
 - \Bigl[ \int_{x^{\prime}}^{2c+x^{\prime}} \hspace{-.1cm} 
             e^{-x^2/2} \frac{dx}{\sqrt{2 \pi}} \Bigr]^2
              \Bigr\} dx^{\prime} 
    :=  (\frac{ \sqrt{N} a}{\sqrt{\pi}} \Bigr) I_{mod}(c)               
\eeq
up to order one errors.  Recall that this holds for either 
$c$ fixed, completing the proof in that case, or $c = c(N) \uparrow \infty$.
In the latter setting, one may easily check that $I_{mod}(c)$ increases
to its value $1 = I_{mod}(\infty)$ which explains the statement at those
scales.

We finish with the case that  $c \downarrow 0$. Here the Edgeworth expansion is used
slightly differently and the growth of the variance is accounted for by
\be
\label{equ:vac}
  V_N(a,c) := \sum_{-N \le \l \le N} \Bigl\{ \int_0^{2c} e^{- \frac{1}{2} (x+ \frac{\l}{\sqrt{N}{a}})^2} 
         \frac{dx}{\sqrt{2 \pi}}     - \Bigl[  \int_0^{2c} e^{- \frac{1}{2} (x+ \frac{\l}{\sqrt{N}{a}})^2}      
                                            \frac{dx}{\sqrt{2 \pi}}   \Bigr]^2   \Bigr\},
\ee
as may be understood upon further review of (\ref{equ:rr}). In this expression it is only the
first sum which plays a role in the asymptotics.  Integrating by parts in that term we find
that
\beqn
\lefteqn{
\sum_{-N \le \l \le N} \int_0^{2c} e^{- \frac{1}{2} (x+ \frac{\l}{\sqrt{N}{a}})^2} 
         \frac{dx}{\sqrt{2 \pi}}} \\
             & = & 
2 c  \sum_{-N \le \l \le N} \frac{e^{- \frac{1}{2}(c+ \frac{\l}{\sqrt{N}{a}})^2}}{\sqrt{2 \pi}} 
                   + O \Bigl( \sqrt{N} c  \int_0^{2c}
                   \sum_{-N \le \l \le N} e^{- \frac{1}{2}(x+ \frac{\l}{\sqrt{N}{a}})^2}  dx  \Bigr) \\
& = & 2 c \sqrt{N} \Bigl[ \int_{-\infty}^{\infty} e^{-x^2/2a^2} \frac{dx}{\sqrt{2 \pi}} \Bigr] (1 + o(1)) + O ( N c^2 ) \\ 
& = & \sqrt{N} (2 a c) (1 + o(1)) = N (b^2 - a^2) (1 + o(1))
\eeqn
after substituting back $ c = \sqrt{N}(b-a) $.  The same computation shows that the 
second, or sum of squares, term in (\ref{equ:vac}) is indeed of lower order. 
The proof is finished.
\bigskip

\noindent{\bf Proof of Theorem \ref{thm:lcov} (radial case)}
Whatever the relation among $a, b, c $ and $d$,
\[
 {\CV}_N \Bigr( \# \I {[a,b]} , \# \I {[c,d]}  \Bigl)
 = \sum_{k=1}^N \Bigl[  p_{N,k} ( [a,b] \cap [c,d] ) - 
                        p_{N,k} ( [a,b] )  p_{N, k} ( [c,d]) \Bigr],
\]
and the estimates obtained in the course of the previous proof may be
revisited. Beginning with the situation $[a^2,b^2] \cap [c^2,d^2] = \varnothing $, 
the covariance reads as  
\be
\label{equ:empty}
 {\CV}_N \Bigr( \# \I{[a,b]} , \# \I {[c,d]} \Bigl)  =   - \sum_{k=1}^N  
         p_{N,k} ( [a^2,b^2] )  p_{N,k} ( [c^2,d^2]).
\ee
Taking for definiteness $b < c$, this sum may be bounded in a simple way as in
\[
 \sum_{k=1}^N   p_{N,k} ( [a^2,b^2] )  p_{N,k} ( [c^2,d^2])  \le 
   \sum_{1 \le  k \le N(c^2+b^2)/2} P (s_k \ge N c^2) + \sum_{N(c^2+b^2)/2 \le k \le N} 
    P(s_k \le N b^2) 
\]
with the first term on the right (being indicative of either) subject to
\[
  \sum_{k=1}^{N(c^2+b^2)/2} P (s_k \ge N c^2) \le  
  e^{-N c^2} \sum_{k =1}^{ N(c^2+b^2)/2} \Bigl(\frac{Nc^2}{k} \Bigr)^k e^k 
  \le c_1 N e^{-N c^2 ( 1 - c^{\prime} + c^{\prime} \log c^{\prime})} . 
\]
with $c^{\prime} = (c^2 + b^2)/2 c^2 < 1$, ensuring
that $1 - c^{\prime} + c^{\prime} \log c^{\prime} > 0$. This second
bound just maximizes the summands, while the first  
is the usual trick of optimizing in Chebychev's inequality 
(as used in (\ref{equ:tail})).  That the right hand side of (\ref{equ:empty}) 
vanishes exponentially fast in $N$ is  established.
The case when $[a^2,b^2] \subset [c^2,d^2]$ follows in kind. Then 
\be
\label{equ:sub}
   {\CV}_N \Bigr( \# \I {[a,b]} , \# \I {[c,d]}  \Bigl) = \sum_{k=1}^N 
                      p_{N,k}([a^2,b^2]) \, (1 -  p_{N,k}([c^2,d^2])), 
\ee
but since now $[a,b]$ and $[c,d]^C$ are disjoint this expression is 
really quite the same as 
(\ref{equ:empty}) and thus
is also exponentially small (though positive) as $N \uparrow \infty$.

The last two cases in which the intervals share a boundary point, either 
$a < b = c < d$ or $a=c$ while say $b < d$, 
also go hand in hand.  In the former we have: with $\d_N = \sqrt{\log N}$,
\beq
\label{equ:over}
\lefteqn{  {\CV}_N \Bigr( \# \I {[a,b]} , \# \I {[b,d]}  \Bigl)
=   - \sum_{k=1}^N p_{N,k}([a^2,b^2]) p_{N,k}([b^2,d^2])} \\ & = &  -
    \sum_{|k - Nb^2| \le \d_N \sqrt{N}} 
    p_{N,k}( [0, b^2] ) p_{N,k}( [b^2, \infty) ) +
    O \Bigl(  \sum_{|k - Nb^2| \ge \d_N \sqrt{N}}   p_{N,k}( [0, b^2] ) p_{N,k}( [b^2, \infty) \Bigr) \nonumber \\ 
& =  & - \frac{b}{2} \sqrt{N}  + o \Bigl( \frac{1}{\sqrt{\log N}} \Bigr) \nonumber.
\eeq
Here the choice of $\d_N$ and the last line are just reprises of the proof of Theorem \ref{thm:num1}.
Finally, if $a = c$ the covariance is 
\[
   {\CV}_N
  \Bigl( \# \I [a,b] , \# \I [a, d] 
  \Bigr) =  \sum_{k=1}^N p_{N,k}([0,a^2]) p_{N,k}([a^2,b^2]) 
  + \sum_{k=1}^N p_{N,k}([a^2,b^2]) p_{N,k}([d^2,\infty]),
\]
where we see that the first sum is of similar type to that in the first line 
of (\ref{equ:over}) and so produces the 
claimed $ \frac{a}{2} \sqrt{N}
+ o(1)$ as $N \uparrow \infty$.  The second term does not affect this outcome;
having the same form as either (\ref{equ:empty}) or (\ref{equ:sub}) it is
exponentially small.  The proof is finished.

\section{Extensions to the Quaternion Ensemble}
\setcounter{num}{5}
\setcounter{equation}{0}

Replacing the 
complex Gaussian entries in the definition of $P_N$ with real quaternion Gaussians produces 
another solvable ensemble $P_N^Q$ first by studied Ginibre.  In this case the eigenvalue density has
the form
\[
 d P_N^Q(z_1, \cdots z_N) = \frac{1}{Z_N^Q} \prod_{1 \le \l < k \le N}
           |z_{\l} - z_k|^2 | z_{\l} - \bar{z_k}|^2 
     \prod_{1 \le k \le N}  |z_k - \bar{z_k} |^2 \, d \tilde{\mu}_N(z_1) \cdots d \tilde{\mu}_N(z_N)  
\]
where now $d {\tilde \mu}(z) = e^{-2 N |z|^2} d \Re(z)d \Im(z)$ is the appropriately scaled weight 
to keep the limiting density of states uniform on $|z| \le 1$.
The introduction of the additional differences-squared terms breaks the rotation symmetry 
of $P_N$. It also breaks some of its analytic structure. While
there again exist explicit formulas for the finite dimensional marginals or correlation functions, these
are now expressed in terms of quaternion determinants or pfaffians.  That is, $P_N^Q$ is not a determinantal point 
field.  Still, part of the previous discussion extends readily
to this model for the reason that the radial components of the eigenvalues under $P_N^Q$ have a nearly
identical statistical description to those in the $P_N$ ensemble. 

The fact is that in either $P_N$ or $P_N^Q$ one may integrate out the angles from the start, producing
a full joint density on the moduli $\{ r_k = |z_k|, k = 1, \dots, N \}$ (see \cite{FS} for related 
remarks). 
The resulting densities possess their own structure, and  for certain analysis it is not 
important whether or not they derived from a determinantal setup.
We start with $P_N$. Returning to the original expression (\ref{equ:def1}), we simply expand
the Vandermonde determinants to find that: with  $\sigma$ and $\psi$ permutations on $N$ letters
and $\t_k = \mbox{arg}(z_k)$,
\beq
\label{equ:comrad}
\lefteqn{d P_N(r_1, \cdots , r_N)} \\ 
& = & \frac{1}{Z_N} \sum_{\sigma, \psi} sgn(\sigma) sgn(\psi) \Bigl\{ \prod_{k=1}^N              
 \int_{-\pi}^{\pi} 
   z_k^{\sigma_{k} -1 } {\bar{z_k}}^{\psi_{k} - 1} d \t_k \Bigr\} \times \prod_{k=1}^N  e^{-N r_k^2} r_k \, d r_k  \nonumber \\
& = &    \frac{1}{N!} \sum_{\sigma}
p_N^{(1)}(r_{\sigma_1}) p_N^{(2)}( r_{\sigma_2}) \cdots p_N^{(N)}
(r_{\sigma_N}) \, d r_1 \cdots d r_N \nonumber
\eeq
in which $p_N^{(\l)}(r) = r^{2 \l -1} e^{-Nr^2}$,
normalized to be a probability density on $r > 0$.
For $P_N^Q$ use is made of the identity
\[
   \det \Bigl[ z_{\l}^{k-1}  \ \ \bar{z_{\l}}^{k-1} \Bigr]_{1 \le k \le 2N, 1 \le \l \le N}
    =  \prod_{1 \le \l < k \le N}
           |z_{\l} - z_k|^2 | z_{\l} - \bar{z_k}|^2 
     \prod_{1 \le k \le N}  (z_k - \bar{z_k} ).
\]
Expanding the left hand side, multiplying by $\prod (z_k - \bar{z_k} ) e^{-N|z_k|^2}$
and integrating over the $\{\t_k \}$ variables 
produces: with $\psi$ now a permutation on $2N$ letters,
\[
d P_N^Q (r_1, \cdots , r_N)
 =  \frac{1}{Z_N^Q} \sum_{\psi} sgn(\psi) \prod_{k=1}^N  \Bigl\{ \int_{-\pi}^{\pi} 
   z_k^{\psi_{2k-1} -1 } {\bar{z_k}}^{\psi_{2k} - 1} (z_k - \bar{z_k}) d \t_k  \Bigr\}
     \times \prod_{k=1}^N e^{-2 N r_k^2} r_k \, d r_k. 
\] 
Noting that only terms with $|\psi_{2k-1} - \psi_{2k}| = 1$ for all $k$ contribute, we continue
the computation as in
\be
\label{equ:quatrad}
d P_N^Q (r_1, \cdots , r_N) 
 =  
\frac{1}{N!} \sum_{ \sigma}
p_{2N}^{(2 \cdot 1)}(r_{\sigma_1}) p_{2N}^{(2 \cdot 2)}( r_{\sigma_2}) \cdots 
p_{2N}^{(2 \cdot N)}
(r_{\sigma_N}) \, d r_1 \cdots d r_N. 
\ee
Here we are back to an $N^{th}$ order permutation in $\sigma$ and the $p_{M}^{(k)}$ have
precisely the same definition as in their appearance in (\ref{equ:comrad}). For background
on this and related computations we refer to \cite{M}, Appendix 18.

Certainly the covariance formula (\ref{equ:cov1}) may be derived directly from
(\ref{equ:comrad}), circumventing the need for computing correlations. It follows from
(\ref{equ:quatrad}) that
the covariances of general linear statistics on the moduli in the $P_N^Q$ ensemble
will satisfy the same formula, so long as each appearance of $\frac{1}{N} s_{\l}$
is replaced with $\frac{1}{2 N} s_{ 2 \l}$.   The asymptotics follow suit
with identical proofs:

\begin{cor} Under $P_N^Q$ the covariances of linear statistics of the radial components
have the identical $N \uparrow \infty$ asymptotic behavior to their counterparts in the
$P_N$ ensemble described in Theorem {\ref{thm:num1}} and the latter half
of Theorem {\ref{thm:lcov}}.
\end{cor}

Further, it should be clear that, by repeating the proof of Theorem {\ref{thm:r1}}, 
we may obtain Gaussian fluctuations for smooth radial linear statistics in $P_N^Q$.  
More interesting, for indicator type test function, while we may no longer 
quote the Costin-Lebowitz result directly, we do have the following.

\begin{thm} {\label{thm:last}}
The properly normalized radial number statistics $\# \I [a,b]$
under $P_N^Q$ also satisfy a central limit theorem as long as $N (b-a) \uparrow \infty$ 
as $N \uparrow \infty$. 
\end{thm}

The proof of Theorem {\ref{thm:last}} is actually a slight reworking of Costin-Lebowitz; 
for that reason it is postponed until the next section.

\bigskip

\noindent{\bf Remark} While given by a {\em permanent}, point processes
of type (\ref{equ:comrad}) or (\ref{equ:quatrad}) 
are not the bosonic fields appearing in the mathematical
physics literature (see \cite{DE1} and references therein).
They are in fact simpler objects; it would be be interesting to know if
something of their kind arose in other contexts.

\section{Appendix}
\setcounter{num}{6}
\setcounter{equation}{0}

For completeness we include the statement and a sketch of the
proof of the
Costin and Lebowitz theorem.  In
keeping with the setting and notation used above, we say 
a measure $P_N$ on $\C^N$ is
a determinantal point field if all $\l$-dimensional 
correlation functions $\rho_{N, {\l}}$ may be written as
\[
   \rho_{N, {\l}} (z_1, z_2, \cdots, z_{\l}) = \frac{N!}{(N- \l)!}
    \det \Bigl[ K_N(z_i, \bar{z_j}) \Bigr]_{1 \le i,j \le \l} \mu_N(dz_1) \cdots \mu_N(dz_{\l})
\]
(compare (\ref{equ:cor})) in which  
$K_N$ is a Hermitian kernel ($K_N(z,w) = \bar{K_N}(w,z)$) 
satisfying $||K_N|| \le 1$ as an operator on $L^2(\C, \mu_N)$
for some measure $\mu_N$.

\begin{thm} {\bf  (O. Costin, J. Lebowitz)} Let
$P_N$ be a determinantal random point field  and let $\# A$ equal the number of
points lying in $A \subset \C$ under that measure. As
long as ${\mbox Var}_N[ \# A] \uparrow \infty$ with $N \uparrow
\infty$ then $ (\# A - E [\# A])/\sqrt{\mbox{Var}_N[\#A]}$
converges in distribution to a Gaussian random variable
of mean zero and variance one as $N \uparrow \infty$.
\end{thm}

\noindent{\bf Proof}  Following \cite{S1} as well as the original
\cite{CL}, the proof 
checks that the cumulants
$C_{\l}^{N}(\# A)$ fall into line as $N$ gets large. Recall
that for any random variable $X$ its $\l$-th cumulant 
$C_{\l}(X)$ is defined
through
\[
  \log E [ e^{i  t X} ] 
   = \sum_{\l = 1}^{\infty}  C_{\l}(X) \frac{ (i t)^{\l}}{\l !}
\]
and that $X$ is Gaussian if and only if $C_{\l}(X) = 0$ for all
$\l \ge 2$ and $C_{2}(X) > 0$ (see \cite{L}). 

Next define $A_N =  K_N {\mathbbm 1}_{A}$, the restriction of $K_N$ to 
$L^2(A, \mu_N)$ and notice that (compare \ref{equ:covfin}) 
$E_N[ \#A] = C_1^N(\# A) = \mbox{Trace} (A_N)$ and ${\mbox Var}_N [ \# A] $
$ = C_2^N(\# A) = \mbox{Trace} (A_N - A_N^2)$.  For the higher
cumulants, there is the following recursion:
\be
\label{equ:cr}
  C_{\l}^N(\#A) = (-1)^{\l} (\l -1)! \, \mbox{Trace} (A_N - A_N^{\l})
                  + \sum_{k =2}^{\l-2} \a_{k \l} C_{k}^N(\# A).
\ee
in which  $\{ \a_{k \l} \}$ are certain combinatorial factors.
The origins of (\ref{equ:cr}) lie in the connection between 
the cumulants and the so-called cluster (or Ursell) functions. The
latter are defined in terms of the correlation functions.
In particular, with  $u_{N,k}$ denoting 
the $k$-point cluster function, 
\be
\label{equ:urs}
  u_{N,k} ( z_1, z_2,\dots, z_k) = \sum (-1)^{k-m} (m - 1)! \, \rho_{N, |S_1|} \cdots  \rho_{N, |S_m|}
\ee
where the sum is over $m \ge 1$ and, for each $m$, all partitions of $\{1,2,\cdots, k\}$ into
(nonempty) subsets $S_1$ etc. (also, $\rho_{|S_1|, N}$ indicated the $|S_1|$-point function
in the variables corresponding to that partition).  Defining also 
\be
\label{equ:burs}
   U_k^N(A) = \int_{A} \cdots \int_{A} u_{N,k}(z_1, z_2, \dots, z_k) \, d x_1 d y_1 \cdots d x_k d y_k, 
\ee  
the advertised connection is provided by the identity
\[
  \sum_{k = 1}^{\infty} C_{k}^N(\# A) \frac{z^k}{k!} = \sum_{k=1}^{\infty} \frac{1}{k!} U_k^N(A) (e^z - 1).
\]
This last display implies that
\be
\label{equ:precr}
 C_{\l}^N(\#A) =
         \sum_{k =2}^{\l-2} \a_{k \l} C_{k}^N(\# A) + (-1)^{\l} (\l - 1)! \, C_1^N(\# A) + U_k^{N}(A).
\ee
And finally, the fact is that, for determinantal point fields,
$U_{\l}^N(A) = (-1)^{\l-1} (\l - 1)! \, \mbox{Trace} ( A_N^{\l} )$, explaining (\ref{equ:cr}).

Now note the simple estimate
\be
\label{equ:ce}
   0 \le \mbox{Trace} (A_N - A_N^{\l}) \le (\l -1) \mbox{Trace}(A_N -A_N^2)
\ee
which may be checked by writing $\mbox{Trace} (A_N - A_N^{\l})$ 
$= \sum_{j=1}^{\l-1}$  $\mbox{Trace} $ $ (A_N^{j} - A_N^{j+1})$ 
$\le  \sum_{j=1}^{\l-1}$ $||A_N^{j-1}|| $ $ \mbox{Trace} (A_N - A_N^{2})$
and recalling that $||K_N || \le 1$. Substituting 
(\ref{equ:ce}) into (\ref{equ:cr}) we have 
that $C_{\l}^N(\# A) = O(C_2^N(\# A))$ for all $\l \ge 2$.  And since
$C_2^N(\# A) \uparrow \infty$ as $N \uparrow \infty$, we also
have that
\[
C_{\l}^N \Bigl( \frac{\# A - E_N[\# A]}{\sqrt{{Var}_N[\# A]}} \Bigr)
 = \frac{C_{\l}^N(\# A)}{(C_2^N(\# A))^{\l/2}} \rightarrow 0
\]
for $\l \ge 2$.  The proof is finished.

\bigskip

\noindent{\bf Proof of Theorem \ref{thm:last}} Consider more generally a set of probability densities
$\{ p_{N,k} \} $ ($k =1,2, \dots N$) on $\C$ and the point process defined for each $N$
by joint density
\[
  P_N(z_1, \dots , z_N) = \frac{1}{N!} \sum_{permutations \, \sigma} p_{N,\sigma(1)} (z_1) 
                                p_{N, \sigma(2)} (z_2)  \cdots p_{N, \sigma(N)} (z_N).
\]
We show that the number statistic for such an ensemble satisfies a central limit theorem
if its variance grows without bound as $N \uparrow \infty$.  This may then by applied to
the particular case of (\ref{equ:quatrad}) in which a choice of $p_{N,k}$ has been made ($ = p_{2N}^{(2k)}$
in that notation), things are
restricted to the positive reals rather than all of $\C$, and the needed growth of the variance 
has already been checked.

This proof also rests on the behavior of the cumulants. While not determinantal,
the correlation functions now have the simpler form:
\[
    \rho_{\l, N}(z_1, \cdots , z_{\l}) = \sum_{1 \le k_1, k_2, \dots k_{\l} \le N} p_{N,k_1} (z_1) p_{N,k_2} (z_2) \cdots 
                                                                                   p_{N,k_{\l}} (z_{\l}).
\]
So, returning to the Costin-Lebowitz proof, we make use not of (\ref{equ:cr}), but instead 
the preliminary cumulant formula (\ref{equ:precr}).
With $\rho_{\l,N}$ as in the previous display, a calculation in (\ref{equ:burs}) 
and (\ref{equ:urs}) shows that 
\[
   U_{\l}^N(A) = (-1)^{\l - 1} (\l - 1)! \, \sum_{k = 1}^{N} [p_{N,k}(A)]^{\l}
\]
in the present case. That is, the growth of 
the $\l^{th}$ cumulant is now tied to $\sum_{k=1}^N ( p_{N,k}(A) - p_{N,k}^{\l}(A) )$.
But, $p - p^{\l} \le (\l-1) (p - p^2)$ for $\l > 2$ and any $p \in [0,1]$, which is to say that $\l^{th}$
cumulant is dominated by a constant multiple of the second. The proof is finished.

\bigskip

\noindent{\bf Acknowledgments} A serious debt of gratitude is due to the referee. 
In addition to helping correct various small glitches, their thoroughness identified 
some important errors in an earlier version of this paper.


\begin{thebibliography}{10}

\bibitem{B}
{\sc Z.D. Bai}:
Circular Law.
{\it Annals of Probability} {\bf 25}, 494-529 (1997).

\bibitem{BS}
{\sc Z.D. Bai, J. Silverstein}:
CLT of linear spectral statistics of large dimensional
sample covariance matrices.
{\it To appear, Annals of Probability}, (2003).

\bibitem{BE} 
{\sc E. Basor}: 
Distribution functions for random variables for ensembles 
of positive Hermitian matrices.
{\it Comm. Mathematical Physics} {\bf 188}, 327-350 (1997).

\bibitem{BW} 
{\sc E. Basor, H. Widom}: 
Determinants of Airy operators and applications to random matrices.
{\it Journal of Statistical Physics} {\bf 96}, 1-20 (1999).


\bibitem{BR}
{\sc R.N. Bhattacharya, R. Ranga Rao}:
Normal Approximation and Asymptotic Expansions.
Wiley, New York, (1976).


\bibitem{BDK}
{\sc D. Bump, P. Diaconis, J. Keller}:
Unitary Correlations and the Fej{\'e}r Kernel.
{\it Mathematical Physics, Analysis, and Geometry} 
{\bf 5}, 101-123 (2002).

\bibitem{CL}
{\sc A. Costin, J.L. Lebowitz}:
Gaussian fluctuations in random matrices.
{\it Physical Review Letters} {\bf 75}, 69-72 (1995).


\bibitem{CT}
{\sc Y.S Chow, H. Teicher}:
{\em Probability Theory: Independence, Interchangeability,
Martingales}  (Second edition), Springer-Verlag (1988).



\bibitem{DE}
{\sc P. Diaconis, S.N. Evans}:
Linear functionals of eigenvalues of random matrices.
{\it Transactions  AMS}
{\bf 353}, 2615-2633 (2001).

\bibitem{DE1}
{\sc P. Diaconis, S.N. Evans}:
Immanants and finite point processes.
{\it J. Comb. Theory Ser. A.} {\bf 91},
305-321 (2000).

\bibitem{E}
{\sc A. Edelman}:
The probability that a random real Gaussian matrix has k real eigenvalues,
related distributions, and the circular law.
{\it Journal of Multivariate Analysis} 
{\bf 60}, 203-232 (1997).


\bibitem{F}
{\sc P.J. Forrester}: 
Fluctuation formula for complex random matrices.
{\it Journal of Physics A: Mathematical and General} 
{\bf 32}, 159-163 (1999).


\bibitem{FS}
{\sc Y.V. Fyodorov, H-J Sommers}:
Random matrices close to Hermitian and unitary:
overview of methods and results.
{\it Journal of Physics A: Mathematical and General}
{\bf 36}, 3303-3347 (2003). 


\bibitem{G}
{\sc J. Ginibre}:
Statistical ensembles of complex, quaternion, and real matrices.
{\it Journal of Mathematical Physics} {\bf 6},  440-449 (1965).

\bibitem{Gu}
{\sc A. Guionnet}:
Large deviations upper bounds and central limit theorems for
 non-commutative functionals of Gaussian large random matrices. 
{\it Ann. Inst. H. Poincar{\'e} Probab. Statist.} {\bf 38},  
341 - 384 (2002).


\bibitem{I}
{\sc S. Israelsson}
Asymptotic fluctuations of a particle system with singular
interaction.
{\it Stoch. Processes and their Appl.} {\bf 93}, 25-56 (2001).





\bibitem{KKP}
{\sc A. M. Khorunzhy, B.A. Khoruzhenko, and L.A. Pastur}:
Asymptotic properties of large random matrices with independent
entries.
{\it Journal of Mathematical Physics} {\bf 37}, 5033-5059 (1996).


\bibitem{Jo1}
{\sc K. Johansson}: 
On the fluctuation of eigenvalues of random Hermitian matrices.
{\it Duke Math Journal }
{\bf 91}, 151-204 (1998).

\bibitem{Jo2}
{\sc K. Johansson}:
On random matrices from the the classical compact 
groups.
{\it Annals of  Mathematics}
{\bf 145}, 519-545 (1997).

\bibitem{KS}
{\sc J.P. Keating, N.C. Snaith}:
Random matrix theory and $\zeta(1/2+it)$.
{\it Communications in Mathematical Physics}
{\bf 214}, 57-89 (2000).

\bibitem{LS}
{\sc N. Lehmann, H.J. Sommers}:
Eigenvalue statistics of random real matrices.
{\it Physics Review Letters} {\bf 67}, 941-944 (1991).

\bibitem{L}
{\sc E. Lukacs}:
{\em Characteristic Functions} (Second edition), Griffin, London (1970).

\bibitem{Ma}
{\sc O. Macchi}:
The coincidence approach to stochastic point processes.
{\it Advances in Applied Probability} {\bf 7}, 83-122 (1975).

\bibitem{M}
{\sc M. L. Mehta}:
{\em Random Matrices} (Second edition), Academic Press, Boston, 1991.


\bibitem{R1}
{\sc B. Rider}:
A limit theorem at the edge of a non-Hermitian random matrix ensemble.
{\it Journal of Physics A: Mathematical and General}
{\bf 36}, 3401 - 3409 (2003).

\bibitem{R2}
{\sc B. Rider}:
Order Statistics and Ginibre's Ensembles.
{\it To appear, Journal of Statistical Physics} (2004).

\bibitem{SiS}
{\sc Ya. Sinai, A. Soshnikov}:
Central limit theorem for traces of
large random symmetric matrices with independent entries.
{\it Bol. Soc. Brasil. Mat.} {\bf 29}, 1-24 (1998).


\bibitem{S1}
{\sc A. Soshnikov}: 
Central limit theorems for local linear 
statistics in classical compact groups and
related combinatorial identities.
{\it Annals of Probability}
{\bf 28}, 1353-1370 (2000).


\bibitem{S2}
{\sc A. Soshnikov}: 
Gaussian fluctuations for Airy, Bessel and
and other determinantal random point fields.
{\it Journal of  Statistical Physics}
{\bf 100}, 491-522 (2000).

\bibitem{S3}
{\sc A. Soshnikov}: 
Gaussian limits for determinantal random point fields. 
{\it Annals of Probability}
{\bf 30}, 171-181 (2002).

\bibitem{S4}
{\sc A. Soshnikov}:
Determinantal random point fields.
{\it Russian Mathematical Surveys}
{\bf 55}, 923-975 (2000).

\bibitem{W}
{\sc K. Wieand}:
Eigenvalue distributions of random unitary matrices.
{\it Probability Theory and Related Fields}
{\bf 123}, 202-224 (2002).


\end{thebibliography}
\end{document}